\documentclass[12pt,A4paper]{article}
\usepackage[english]{babel}
\usepackage{float} 
\usepackage{array}
\usepackage[utf8]{inputenc}
\usepackage{amssymb,amsfonts,amsthm,amsmath}
\usepackage{enumerate}
\usepackage{tikz}
\usepackage{tikz-cd}
\usepackage{anysize}
\usepackage{indentfirst}
\usepackage{xypic}
\usepackage{mathabx}
\usepackage{xy}

\theoremstyle{theorem}
\newtheorem{theorem}{Theorem}

\newtheorem{theo}{Theorem}[section]
\newtheorem{prop}[theo]{Proposition}
\newtheorem{lemma}[theo]{Lemma}

\newtheorem{cor}[theo]{Corollary}
\theoremstyle{definition}
\newtheorem{df}[theo]{Definition}

\newtheorem{obs}[theo]{Remark}

\newcommand{\G}{\Gamma}
\newcommand{\D}{\Delta}
\newcommand{\ra}{\rightarrow}
\newcommand{\f}{\varphi}
\newcommand{\li}{\varprojlim}

\newcommand{\ZC}{Z_{\hat{\mathcal{C}}}}
\newcommand{\CA}{\mathcal{C}}

\newcommand{\bs}{\backslash}
\newcommand{\GA}{\mathcal{G}}
\newcommand{\HA}{\mathcal{H}}
\newcommand{\U}{\mathcal{U}}
\newcommand{\FG}{\Pi_1^{\CA}(\GA,\G)}

\newcommand{\FGP}{\pi_1^{\CA}(\G)}
\newcommand{\GG}{(\GA,\pi,\G)}
\newcommand{\SG}{S^{\CA}(\GA,\G,H)}

\def\dbigcup{\mathinner{\bigcup \mkern -13.2mu \rlap{\raise 0.6ex\hbox{.}}\mkern 14.9mu}}

\title{The profinite completion of the fundamental group of infinite graphs of groups}
\author{Mattheus P. S. Aguiar and Pavel A. Zalesskii}

\begin{document}
	\pretolerance10000
	\maketitle
	
	\begin{abstract}
		Given an infinite graph of profinite groups $(\GA,\G)$ we construct a profinite graph of groups $(\overline{\GA},\overline{\G})$ such that $\G$ is densely embedded in $\overline{\G}$, the fundamental profinite group $\Pi_1(\overline{\GA},\overline{\G})$ is the profinite completion of $\pi_1(\GA,\G)$ and the standard tree $S(\GA,\G)$ embeds densely in the standard profinite tree $S(\overline{\GA},\overline{\G})$. This answers a Ribes' question  \cite[Question 6.7.1]{Ribes}. Generalizing the main results of  \cite{Zalesskii} and \cite{Chagas} we answer two other questions of Ribes  \cite[Questions 15.11.10 and 15.11.11]{Ribes}   proving that a virtually free group $G$  is subgroup conjugacy separable and the normalizer $N_G(H)$ of a finitely generated subgroup  $H$  of $G$ is dense in  $N_{\widehat{G}}(\overline{H})$.  We also give an entirely new description of the fundamental group of a profinite graph of groups using the language of paths.
	\end{abstract}
	
	
	\section{Introduction} 

	The  Bass-Serre theory of groups acting on trees  appeared   in 1968-1971 and changed completely  the face of combinatorial group theory. It gave understanding that the main constructions of classical combinatorial group theory can be expressed uniformly as groups acting on trees and it allowed to use the metric of the tree instead  of manipulation with words that represent elements of the group.
	
	Approaching 80's  the need of such theory appeared in the category of profinite groups, due to the   difficulty of obtaining basic results on amalgamated free products or HNN-extensions by existing methods. An element of a profinite group can not be written as a word of generators  and so the classical methods of combinatorial group theory are absent in the  class of profinite groups. The profinite version of the Bass-Serre theory was initiated in \cite{GR} and then developed in a series of papers by Melnikov, Ribes and the second author. In 2017 was published a book \cite{Ribes} on the subject where the reader can find basics and the main results of the theory. The third part of the book is dedicated to applications to abstract groups, where the methods are based on an interplay between the classical Bass-Serre theory and its profinite version.   More concretely, if $G$ is the fundamental group $G=\pi_1(\GA,\Gamma)$ of a finite graph of groups $(\GA, \Gamma)$ then the profinite completion $\widehat G$ is the fundamental group $\widehat G=\Pi_1(\widehat \GA,\Gamma)$ over the same graph with fibers being the profinite completions. However, it works only when the underlying graph $\Gamma$ is finite, and sometimes, for example when $G$ is  infinitely generated free-by-finite, the underlying graph  $\Gamma$ is infinite.
	
	This motivated L. Ribes to ask weather there is a reasonable way to define a graph of groups $(\GA', \overline\Gamma)$ such that $\Pi_1(\GA',\overline\Gamma)=\widehat{\pi_1(\GA,\Gamma)}$ when $\Gamma$ is infinite and the vertex groups are finite (see Open Question 6.7.1 \cite{Ribes} for the precise formulation).
	
	Our first result answers this question; in fact we state it more generally, namely for the pro-$\CA$ completion, where $\CA$ is the class of finite groups closed for subgroups and extensions.  
	
	\begin{theorem} \label{B}
		Let $(\GA,\G)$ be an abstract reduced graph of finite groups over an abstract graph $\G$ such that $\pi_1(\GA,\G)$ is residually $\CA$. Then there exists a profinite graph $\overline{\G}$ that contains $\G$ as a dense subgraph and a profinite graph of finite groups $(\overline{\GA},\overline{\G})$, with $\overline{\GA}(m)=\GA(m)$, whenever $m \in \G$, such that
		\begin{enumerate}[(a)]
			\item The pro-$\CA$ fundamental group $\Pi_1^{\CA}(\overline{\GA},\overline{\G})$ is the pro-$\CA$ completion of $\pi_1(\GA,\G)$;
			\item The standard tree $S^{abs}(\GA,\G)$ for $G$  is densely embedded in the standard pro-$\CA$ tree $S^{\CA}(\overline{\GA},\overline{\G})$ in such a way that the action of $\Pi_1^{\CA}(\overline{\GA},\overline{\G})$ on $S^{\CA}(\overline{\GA},\overline{\G})$ extends the natural action of $\pi_1(\GA,\G)$ on $S^{abs}(\GA,\G)$. 
		\end{enumerate}
	\end{theorem}	
	
	This theorem allows us to answer Open Question $15.11.10$ of \cite{Ribes} on the density of the normaliser of a finitely generated subgroup of a virtually free group in the normaliser of its closure in the profinite completion. 
	
	\begin{theorem} \label{A}
		Let $G$ be a virtually free abstract group and $H$ a finitely generated subgroup of $G$.  Then \[\overline{N_G(H)}=N_{\widehat G}(\overline{H}).\]
	\end{theorem}
	
	Theorem \ref{A} generalises the main result of \cite{Zalesskii}, where it was proved for the finitely generated case. In fact we prove a more general  result, namely  for the pro-$\CA$ completion of a free-by-$\CA$ group assuming that $H$ is closed in the pro-$\CA$ topology (see Theorem \ref{a}).
	
	\medskip	 
	
	Let $G=\pi_1(\GA,\Gamma)$ be the fundamental group of  a graph of groups and suppose $G$ is residually $\CA$. Let $\overline{\mathcal{G}}(m)$  be the closure of $\GA(m)$ in $G_{\widehat\CA}$.
	We prove a generalisation of Theorem \ref{B} by showing that we can define a profinite graph $\overline{\G}$ that contains $\G$ as a dense subgraph and a profinite graph of groups $(\overline{\GA},\overline{\G})$, where the groups $\GA(m)$ are not necessarily finite. 
	
	\begin{theorem} \label{C}
		Let $(\GA,\G)$ be a reduced graph of  groups   and $G=\pi_1(\GA,\G)$ be its fundamental group. Assume that $G$ is residually $\CA$.   Then,
		\begin{enumerate}[(a)]
			\item There exists a profinite graph of pro-$\CA$ groups $(\overline{\GA},\overline{\G})$ such that $\G$ is densely embedded in $\overline{\G}$;
			\item for each $m\in\Gamma$ the vertex group is $\overline\GA(m)$;
			\item The fundamental pro-$\CA$ group $\Pi=\Pi_1^{\CA}(\overline{\GA},\overline{\G})$ of $(\overline{\GA},\overline{\G})$ is the pro-$\CA$ completion of $G$,  so that   all the vertex groups of   $(\overline{\GA},\overline{\G})$ embed in  $\Pi_1^{\CA}(\overline{\GA},\overline{\G})$ ( i.e. $(\overline{\GA},\overline{\G})$  is injective in the terminology of \cite{Ribes}).
			\item If in addition we assume that $\GA(m)$ is closed in the pro-$\CA$ topology of $G$, for every $m \in \G$, the standard tree $S^{abs}=S(\GA,\G)$ of the graph of groups $(\GA,\G)$  embeds densely  in the standard pro-$\CA$-tree $S=S^{\CA}(\overline{\GA},\overline{\G})$ of the profinite graph of profinite groups $(\overline{\GA},\overline{\G})$ in such a way that the action of $\Pi_1^{\CA}(\overline{\GA},\overline{\G})$ on $S^{\CA}(\overline{\GA},\overline{\G})$ extends the natural action of $\pi_1(\GA,\G)$ on $S^{abs}(\GA,\G)$. 
		\end{enumerate}
	\end{theorem}  
	
	In order to prove this theorem we give a new description of the profinite fundamental group of a profinite graph of profinite groups; namely we  transport to the profinite context the definition of the fundamental group of a graph of groups in the language of paths. It gives the possibility to work with morphisms of profinite graphs of groups and to use their projective limits.
	
	\smallskip
	We  apply  the above results to establish subgroup conjugacy separability of virtually free groups answering Open Question 15.11.11 of \cite{Ribes}.  A group $G$ is said to be subgroup conjugacy separable if  whenever $H_1$ and $H_2$ are finitely generated  subgroups of $G$  then $H_1$ and $H_2$ are conjugate in $G$ if and only if their images in every finite quotient  are conjugate. 
	
	\begin{theorem} \label{D}  
		Let $G$ be a virtually free  group. Then $G$ is subgroup conjugacy separable.
	\end{theorem}
	
	In fact we prove the subgroup conjugacy $\CA$-separability of a residually $\CA$ free-by-$\CA$ abstract group $G$. This generalises the main result of \cite{Chagas} where it is proved for finitely generated free-by-$\CA$ abstract groups. 
	
	This paper is organised as follows. In Section $2$ we  recall the elements of the Bass-Serre theory for abstract and profinite graphs that will be used through the text and establish the notation to be used in subsequent sections. In Section 3 we introduce a new definition of the fundamental group of a graph of pro-$\CA$-groups at based point and prove that it is equivalent to the one that appears in \cite{Ribes}. Section $4$ is the technical heart of the paper and contains the proof of Theorem \ref{C}. We construct explicitly the graph of profinite groups $(\overline \GA, \overline\Gamma)$ starting with an infinite graph of groups $(\GA,\G)$. With this construction in hand, Theorem \ref{B} follows immediately. In Section $5$ we prove Theorem \ref{A} on the closure of normalisers. The techniques here are largely based on the abstract and profinite Bass-Serre theory and their close interrelation. We finish the paper with Section $6$ that contains the proof of Theorem \ref{D}.   	
	
	\section*{Notation}
	
	In this paper we shall follow the notation of \cite{RZ} and \cite{Ribes}, where the reader can find the main concepts and results related to profinite graphs and groups. Throughout the article one assumes that $\CA$ is a pseudovariety of finite groups closed for extensions, i.e., $\CA$ is a non-empty collection of (isomorphism classes of) finite groups closed for taking subgroups, homomorphic images and extensions. For example, $\CA$ can be the class of all finite groups, the class of all finite p-groups for a fixed prime number $p$, or the class of all finite solvable groups. A 'pro-$\CA$' group $G$ is an inverse limit of groups in $\CA$; this is a compact, Hausdorff and totally disconnected topological group with the property $G/U \in \CA$, whenever $U$ is an open normal subgroup of $G$.   
	
	Let $G$ be an abstract group. Recall that the (full) 'pro-$\CA$ topology' of $G$ is the unique topology that makes $G$ into a topological group in such a way that the set $\mathcal{U}$ of all normal subgroups $U$ of $G$ with $G/U \in \CA$ form a fundamental system of neighbourhoods of the identity element 1. One says that $G$ is 'residually $\CA$' if this topology is Hausdorff, i.e., if $\bigcap_{U \in \mathcal{U}} U=1$. The 'pro-$\CA$ completion' $G_{\widehat{\CA}}$ of $G$ is the pro-$\CA$ group \[G_{\widehat{\CA}}=\li_{U \in \mathcal{U}}G/U.\] The natural homomorphism $G \ra G_{\widehat{\CA}}$ is continuous ($G$ is endowed with its pro-$\CA$ topology). If $G$ is residually $\CA$, this homomorphism is an injection, and one identifies $G$ with its image in $G_{\widehat{\CA}}$, so that $G \leq G_{\widehat{\CA}}$. In this case the topology on $G$ induced by the topology of $G_{\widehat{\CA}}$ is precisely its full pro-$\CA$ topology. If $X \subseteq G$, one denotes the closure of $X$ in $G$ by $Cl(X)$, and the closure of $X$ in $G_{\widehat{\CA}}$ by $\overline{X}$; one notes that  $Cl(X)= G \cap \overline{X}$ and $\overline{X}=\overline{Cl(X)}$, (cf. \cite[Section 3]{RZ94}).
	
	One shall often be interested in abstract groups $G$ that are free-by-$\CA$; that is, $G$ contains a normal free subgroup $\Phi$ such that $G/\Phi \in \CA$, or equivalently, $\Phi$ is open in the pro-$\CA$ topology of $G$. When $\CA$ is the class of all finite groups, one reverts to the usual terminology 'free-by-finite' (or 'virtually free') rather than free-by-$\CA$.
	
	\section{Preliminaries}
	
	In this section we recall the necessary notions of the Bass-Serre theory for abstract and profinite graphs.
	
	\begin{df}[(Profinite) graph]
		A (profinite) graph is a (profinite space) set $\Gamma$ with a distinguished closed nonempty subset $V(\Gamma)$ called the vertex set, $E(\Gamma)=\Gamma-V(\Gamma)$ the edge set and two (continuous) maps $d_0,d_1:\Gamma \rightarrow V(\Gamma)$ whose restrictions to $V(\Gamma)$ are the identity map $id_{V(\Gamma)}$. We refer to $d_0$ and $d_1$ as the incidence maps of the (profinite) graph $\Gamma$.  
	\end{df}
	
	A morphism $\alpha:\G \ra \D$ of profinite graphs is a continuous map with $\alpha d_i=d_i \alpha$ for $i=0,1$. By \cite[Proposition 2.1.4]{Ribes} every profinite graph $\G$ is an inverse limit of finite quotient graphs of $\G$.
	
	\begin{df}[pro-$\CA$-tree] Let $\G$ be a profinite graph. Define $E^{*}(\G)=\G/V(\G)$ to be the quotient space of $\G$ (viewed as a profinite space) modulo the subspace of vertices $V(\G)$. Let $R$ be a profinite ring and consider the free profinite $R$-modules $[[R(E^{*}(\G),*)]]$ and $[[RV(\G)]]$ on the pointed profinite space $(E^{*}(\G),*)$ and on the profinite space $V(\G)$, respectively. Denote by $C(\G,R)$ the chain complex 
		\begin{equation*}
			\begin{tikzcd}
				0 \arrow{r} & \left[\left[R(E^{*}(\G),*)\right]\right] \arrow{r}{d} & \left[\left[RV(\G)\right]\right] \arrow{r}{\varepsilon} & R \arrow{r} & 0
			\end{tikzcd}
		\end{equation*}
		of free profinite $R$-modules and continuous $R$-homomorphisms $d$ and $\varepsilon$ determined by $\varepsilon(v)=1$, for every $v \in V(\G)$, $d(\overline{e})=d_1(e)-d_0(e)$, where $\overline{e}$ is the image of an edge $e \in E(\G)$ in the quotient space $E^{*}(\G)$, and $d(*)=0$. Obviously $Im(d) \subseteq Ker(\varepsilon)$ and one defines the homology groups of $\G$ as the homology groups of the chain complex $C(\G,R)$ in the usual way:
		\begin{center}
			$H_0(\G,R)=Ker(\varepsilon)/Im(d)$ and $H_1(\G,R)=Ker(d)$.
		\end{center}
		One says that $\G$ is a pro-$\CA$ tree if the sequence $C(\G,\ZC)$ is exact. 
	\end{df}	
	
	By \cite[Proposition 2.3.2]{Ribes} a profinite graph is connected if and only if $H_0(\G,R)=0$, independently of the choice of the profinite ring $R$. Therefore a profinite graph $\G$ is a pro-$\CA$ tree if and only if it is connected and $H_1(\G,\ZC)=0$.  If $v$ and $w$ are elements of a tree (respectively $\CA$-tree) $T$, one denotes by $[v,w]$ the smallest subtree (respectively $\CA$-tree) of $T$ containing $v$ and $w$.
	
	\begin{df}[Galois covering]
		Let $G$ be a pro-$\CA$ group that acts freely on a profinite graph $\Gamma$. The natural epimorphism of profinite graphs $\zeta: \Gamma \rightarrow \Delta=G \backslash \Gamma$ of $\Gamma$ onto the quotient graph by the action of $G$, $\Delta=G \backslash \Gamma$ is called a Galois $\CA$-covering of the profinite graph $\Delta$. The associated pro-$\CA$ group $G$ is called the group associated with the Galois covering $\zeta$ and we denote it by $G=G(\Gamma | \Delta)$. If $\G$ is finite, one says that the Galois covering $\zeta$ is finite. The Galois covering is said to be connected if $\G$ is connected.  If $\Gamma$ does not have non-trivial Galois $\CA$-coverings, then $\zeta$ is called {\it universal} $\CA$-covering and in this case $G=\pi_1^{\CA}(\Gamma)$ is the pro-$\CA$ fundamental group of $\Delta$ (see \cite[Theorem 3.7.1]{Ribes} or  \cite[Theorem 2.8]{Zalesskii}).
	\end{df}

	\begin{df}[Sheaf of (profinite) groups]
		Let $T$ be a (profinite space) set. A sheaf of (profinite) groups over $T$ is a triple $(\GA,\pi,T)$, where $\GA$ is a (profinite space) set and $\pi:\GA \ra T$ is a (continuous) surjection satisfying the following conditions:
		\begin{enumerate}[(a)]
			\item For every $t \in T$, the fiber $\GA(t)=\pi^{-1}(t)$ over $t$ is a (profinite) group (whose topology is induced by the topology of $\GA$ as the subspace topology);
			
			\begin{center}
				and only for the profinite case,
			\end{center}
			\item If we define \[\GA^2=\{(g,h) \in \GA \times \GA \mid \pi(g)=\pi(h)\},\] then the map $\mu:\GA^2 \ra \GA$ given by $\mu_{\GA}(g,h)=gh^{-1}$ is continuous.
		\end{enumerate} 
	\end{df}
	
	\begin{df}
		A morphism $\underline{\alpha}=(\alpha,\alpha'):(\GA,\pi,T) \ra (\GA',\pi',T')$ of sheaves of (profinite) groups consists of a pair of (continuous) maps $\alpha: \GA \ra \GA'$ and $\alpha': T \ra T'$ such that the diagram
		\begin{equation*}
			\begin{tikzcd}
				\GA \arrow{r}{\alpha} \arrow{d}[swap]{\pi} & \GA' \arrow{d}{\pi'}\\
				T \arrow{r}[swap]{\alpha'} & T'  
			\end{tikzcd}
		\end{equation*}
		commutes and the restriction of $\alpha$ to $\GA(t)$ is a homomorphism from $\GA(t)$ into $\GA'(\alpha'(t))$, for each $t \in T$.
	\end{df}	
	
	If $\alpha$ and $\alpha'$ are injective, the morphism $\underline{\alpha}$ is said to be a monomorphism and the image of a monomorphism $\GA \ra \GA'$ is called a subsheaf of the sheaf $\GA'$. If $(\GA,\pi,T)$ is a sheaf and $T'$ is a (closed) subspace of $T$, then the triple \[(\pi^{-1}(T'),\pi|_{\pi^{-1}(T')},T')\] is a subsheaf of $(\GA,\pi,T)$.

	\begin{df}[(Profinite) graph of (profinite) groups]
		Let $\G$ be a connected (profinite) graph with incidence maps $d_0,d_1: \G \ra V(\G)$. A (profinite) graph of groups over $\G$ is a sheaf $(\GA,\pi,\G)$ of (profinite) groups over $\G$ together with two morphisms of sheaves $(\partial_i,d_i):(\GA,\pi,\G) \ra (\GA_V,\pi,V(\G))$, where $(\GA_V,\pi,V(\G))$ is a restriction sheaf of $(\GA,\pi,\G)$ and the restriction of $\partial_i$ to $\GA_v$ is the identity map $id_{\GA_V}$, $i=0,1$; in addition, we assume that the restriction of $\partial_i$ to each fiber $\GA(m)$ is an injection.
	\end{df}
	
	\begin{df} A morphism of graphs of groups $\underline\nu=(\nu,\nu'):(\GA,\Gamma) \longrightarrow (\HA, \Delta)$ is a morphism of sheaves such that  $\nu\partial_i=\partial_i\nu$.
	\end{df}

	Let $I$ be a partially ordered set, $\{(\GA_i,\pi_i,\G_i),\nu_{ij}\}$ an inverse system of finite  graphs of $\CA$-groups. Then $(\GA,\pi,\G)=\li_{i \in I}(\GA_i,\pi_i,\G_i)$ is a profinite graph of pro-$\CA$ groups.

	\medskip
	We shall need later the definition of the graph of pro-$\CA$ groups $(\widetilde \GA, \widetilde\Gamma)$ associated to a profinite graph of pro-$\CA$-groups $(\GA, \Gamma)$ that appears in \cite[Subsection 2.5]{ZM} and \cite[page 185]{Ribes}.
	
	\begin{df}\label{universal graph of groups}
		Let $(\GA, \Gamma)$ be a profinite graph of pro-$\CA$ groups. Let $\widetilde{\G}$ be the universal $\CA$-covering of the profinite graph $\G$ and $\zeta: \widetilde{\G} \ra \G$ be the covering map. Consider the pull-back 
		\begin{equation} \label{pb2}
			\begin{tikzcd}
				\widetilde{\GA} \arrow{d}[swap]{\tilde{\zeta}} \arrow{r}{\tilde{\pi}} & \widetilde{\G} \arrow{d}{\zeta}\\
				\GA \arrow{r}[swap]{\pi} & \G  
			\end{tikzcd}
		\end{equation}
		of the maps $\pi: \GA \ra \G$ and $\zeta:\widetilde{\G} \ra \G$. This means that $\widetilde{\GA}$ is defined by  \[\widetilde{\GA}=\{(x,\tilde{m}) \in \GA \times \widetilde{\G} \mid \pi(x)=\zeta(\tilde{m}),x \in \GA, \tilde{m} \in \widetilde{\G}\} \subseteq \GA \times \widetilde{\G}.\]
		
		Let $\tilde{\pi}:\widetilde{\GA} \ra \widetilde{\G}$ and $\tilde{\zeta}:\widetilde{\GA} \ra \GA$  be the restrictions to $\widetilde{\GA}$ of the canonical projections from $\GA \times \widetilde{\G}$ to $\widetilde{\G}$ and $\GA$ respectively. For a fixed $\tilde{m} \in \widetilde{\G}$, define \[\widetilde{\GA}(\tilde{m})=\pi^{-1}(\tilde{m})=\GA(\zeta(\tilde{m})) \times \{\tilde{m}\},\] which is a group with respect to the operation $(x,\tilde{m})(x',\tilde{m})=(xx',\tilde{m})$, $x,x' \in \GA(\zeta(\tilde{m}))$ and clearly $\widetilde{\GA}(\tilde{m}) \cong \GA(\zeta(\tilde{m}))$.
		
		Then $(\widetilde{\GA},\widetilde{\G})$ is a profinite graph of the  above defined groups $\widetilde{\GA}(\tilde{m})$ over the universal covering graph $\widetilde{\G}$ of $\G$ with boudary maps  $\partial_i(g)=(\partial_i(\tilde\zeta(g)), d_i(\tilde m))$ for $g\in \widetilde\GA(\tilde m)$.
		
	\end{df}
	
	\begin{df}[$G$-transversal] \label{trans}
		Let $G$ be a profinite group that acts on a connected profinite graph $\G$, and let $\f:\G \ra \D=G\backslash\G$ be the canonical quotient map. A $G$-transversal or a transversal of $\f$ is a closed subset $J$ of $\G$ such that $\f_{|J}:J \ra \D$ is a homeomorphism. Associated with this transversal there is a continuous $G$-section or section of $\f$, $j:\D \ra \G$, i.e., a continuous mapping such that $\f j=id_{\D}$ and $j(\D)=J$. 
	\end{df}
	
	Note that, in general, $J$ is not a graph.
	
	\begin{df}[0-transversal]
		We say that a transversal $J$ is a 0-transversal if $d_0(m) \in J$, for each $m \in J$; in this case we refer to $j$ as a 0-section.
	\end{df}
	
	Note that if $j$ is a 0-section, then $jd_0=d_0j$.	
	
	\begin{df}[$J$-specialisation] \label{spec}
		Given a pro-$\CA$ group $H$, define a $J$-specialisation of the graph of pro-$\CA$ groups $(\GA,\pi,\G)$ in $H$ to consist of a pair $(\beta,\beta')$, where $\beta:(\GA,\pi,\G) \ra H$ is a morphism from the sheaf $(\GA,\pi,\G)$ to $H$, and where $\beta':\pi_1^{\CA}(\G) \ra H$ is a continuous homomorphism satisfying the following conditions: 
		\begin{equation*} \label{specialisation}
			\beta(x)=\beta\partial_0(x)=(\beta'\chi(m))(\beta\partial_1(x))(\beta'\chi(m))^{-1}
		\end{equation*}
		for all $x \in \GA$. Here $m=\pi(x)$ and the continuous map $\chi: \G \ra \FGP$ is defined in the following way: $\chi(m)$ is the unique element of $\pi_1^{\CA}(\G)$ such that $\chi(m)(jd_1(m))=d_1j(m).$
	\end{df}
	
	The definition of the pro-$\CA$ fundamental group of a profinite graph of pro-$\CA$ groups used in the literature until present paper is as follows.  	
	
	\begin{df}[The fundamental group of a graph of pro-$\CA$ groups] \label{fgg}
		Choose a continuous 0-section $j$ of the universal Galois $\CA$-covering $\zeta:\widetilde{\G} \ra \G$ and denote by $J=j(\G)$ the corresponding 0-transversal. We define a fundamental pro-$\CA$ group of the graph of groups $(\GA,\pi,\G)$ with respect to the 0-transversal $J$ to be a pro-$\CA$ group $\Pi_1^{\CA}(\GA,\G)$, together with a $J$-specialisation $(\nu,\nu')$ of $(\GA,\pi,\G)$ in $\FG$ satisfying the following universal property:
		
		\begin{center}
			\begin{tikzpicture}
				\node (0) at (0,0) {$H$};
				\node (1) at (-2,2) {$\GA$};
				\node (2) at (2,2) {$\FGP$};
				\node (3) at (0,4) {$\FG$};
				
				\draw[->] (3) edge[dotted] node[left] {$\delta$} (0);
				\draw[->] (1) edge node[left] {$\beta$} (0);
				\draw[->] (2) edge node[right] {$\beta'$} (0);
				\draw[->] (1) edge node[left] {$\nu$} (3);
				\draw[->] (2) edge node[right] {$\nu'$} (3);
			\end{tikzpicture}
		\end{center}
		whenever $H$ is a pro-$\CA$ group and $(\beta,\beta')$ a $J$-specialisation of $(\GA,\pi,\G)$ in $H$, there exists a unique continuous homomorphism $$\delta:\FG \ra H$$ such that $\delta\nu=\beta$ and $\delta\nu'=\beta'$. We refer to $(\nu,\nu')$ as a universal $J$-specialisation of $(\GA,\pi,\G)$. 
	\end{df}
	
	Note that in the pro-$\CA$ case the vertex groups do not necessarily embed in the fundamental group (namely, $\nu$ restricted on fibers might not be a monomorphism); when it is the case, the profinite graph of pro-$\CA$ groups is called {\it injective}. Observe that a profinite graph of pro-$\CA$ group $(\GA,\Gamma)$ can   be replaced with a natural graph of quotient groups $(\overline\GA, \Gamma)$ by replacing $\GA(m)$ by its image $\nu(\GA(m)$ in $\Pi_1^{\CA}(\GA,\Gamma)$ for every $m\in\Gamma$; then  $\Pi_1^{\CA}(\GA,\Gamma)=\Pi_1^{\CA}(\overline\GA,\Gamma)$ and $(\overline\GA,\Gamma)$ is injective (see \cite[Section 6.4]{Ribes} for details). This means that we do not loose generality restricting our attention to injective profinite graphs of pro-$\CA$ groups. 
	
	By \cite[Theorem 6.2.4]{Ribes} the definition of $\FG$ does not depend on the choice of the transversal $J$.
	
	We next show that an injective profinite graph of pro-$\CA$ groups decomposes as an inverse limite of finite graphs of finite $\CA$-groups. For this we shall need the following lemmas.
	
	\begin{lemma}\label{factors through}  Let $(\GA,\G)$ be a profinite graph of pro-$\CA$-groups. Suppose there exists a decomposition  $(\GA,\G)=\varprojlim_{i\in I} (\GA_i,\G_i)$ as a surjective inverse limit of finite graphs of finite $\CA$-groups. Let $\underline\alpha=(\alpha,\alpha'):(\GA, \Gamma)\longrightarrow (\HA, \Delta)$ be a surjective morphism to a finite graph of finite $\CA$-groups and $\eta:(\GA, \Gamma)\longrightarrow H$ be a fiber homomorphism to a finite $\CA$-group $H$. Then $\underline\alpha, \eta$ factor via some $(\GA_k,\G_k)$, i.e. there exist a morphism $\underline\alpha_k:(\GA_k,\G_k)\longrightarrow (\HA, \Delta)$ and a fiber homomorphism $\eta_k:(\GA_k, \Gamma_k)\longrightarrow H$ such that $\underline\alpha=\underline\alpha_k\underline\pi_k$ and $\eta=\eta_k\underline\pi_k$, where $\underline\pi_i=(\pi_i,\pi'_i):(\GA,\G)\longrightarrow (\GA_i,\G_i)$ is the natural projection.\end{lemma}   
	
	\begin{proof} We adapt the proof of  \cite[Lemma 2.1.5]{Ribes}). 
		
		Let $S$ be the equivalence relation on $(\GA,\Gamma)$ whose equivalence classes  are the clopen sets $ \alpha^{-1}(h), h \in \HA$ and $R$ the equivalence relation on $(\GA,\Gamma)$ whose equivalence classes  are the clopen sets $\alpha'^{-1}(m), m \in \Delta$; then $(\GA/S,\Gamma/R)=(\HA, \Delta)$ and $\underline\alpha=(\alpha,\alpha')$ is the natural projection $(\GA, \Gamma)\longrightarrow (\GA/S,\Gamma/R)$. Similarly, for $i\in I$, let $S_i$ be the equivalence relation on $(\GA,\Gamma)$ whose equivalence classes  are the clopen sets $ \pi_i^{-1}(h), h \in \GA_i$ and $R_i$ the equivalence relation on $(\GA,\Gamma)$ whose equivalence classes  are the clopen sets $(\pi'_{i})^{-1}(m), m \in \Gamma_i$, so that $\underline\pi_i=(\pi_i,\pi'_{i}):(\GA,\G)\longrightarrow (\GA/S_i,\Gamma/R_i)$ is the natural projection. Since $(\GA,\G)=\varprojlim_{i\in I} (\GA_i,\G_i)$, the intersections  $\bigcap_{i\in I} S_i$, $\bigcap_{i\in I} R_i$ are the diagonal subsets of $\GA\times \GA$ and $\Gamma \times \Gamma$, respectively. Note that $S,S_i$ are clopen subsets of  $\GA\times \GA$ and $R,R_i$ are clopen subsets of $\Gamma \times \Gamma$. Hence, it follows from compactness of  $\GA\times \GA$ and  $\G\times \G$ that there exists a finite subset $J$ of $I$ such that   $\bigcap_{j\in J} S_j \subseteq S$ and $\bigcap_{j\in J} R_j \subseteq R$. Since poset $I$ is directed, there exists a $k\in I$ such that $S_k \subseteq \bigcap_{j\in J} S_j \subseteq S$ and $R_k\subseteq \bigcap_{j\in J} R_j \subseteq R$. This means that there exists a morphism of graphs of groups $\underline\alpha_k:(\GA_k,\G_k)\longrightarrow (\HA, \Delta)$ such that $\underline\alpha=\underline\alpha_k\underline\pi_k$. 
		
		The statement about a fiber homomorphism $\eta_k:(\GA_k, \Gamma_k)\longrightarrow H$ follows from the first, since we can consider $H$ as a graph of groups $(H,\{v\})$ with underlying graph being one vertex $v$.
	\end{proof}
	
	\begin{lemma}\label{inverse system} Let $(\GA,\G)$ be a profinite graph of pro-$\CA$-groups. Suppose there exists a decomposition  $(\GA,\G)=\varprojlim_{U\in \U} (\GA_U,\G)$ as an inverse limit of profinite graph of finite quotient $\CA$-groups over the same graph $\Gamma$. Assume further that for each $U\in \U$ the graph of groups $(\GA_U,\G)$ decomposes as a surjective inverse limit $(\GA_U,\G)=\varprojlim_{i\in  I_U} (\GA_{i,U},\Gamma_i)$. Then $(\GA_{i,U},\Gamma_i)$ form naturaly an inverse system such that $(\GA,\Gamma)=\varprojlim_{U \in \U, i\in  I_U} (\GA_{i,U},\Gamma_i)$.
		
	\end{lemma}

	\begin{proof} We follow the proof of \cite[Propisition 3.1.3]{Ribes} making all the appropriate changes. Denote by $\varphi_{i,U}:(\GA_U, \Gamma)\longrightarrow (\GA_{i,U},\Gamma_i)$ the canonical projection. Define the indexing set $I=\bigcup_{U\in \U} I_U$. We relabel the elements of $I_U$: an element $i\in I_U$ will be denoted from now on by $(i,U)$.  If $(i,U),(j,V) \in I$, we say $(i,U)\leq (j,V)$ if $U\leq V$ and there exists a morphism of graph of groups 
		$\alpha:(\GA_{j,V},\Gamma_j)\longrightarrow  (\GA_{i,U}, \Gamma_i)$ such that the diagram 
		
		$$\xymatrix{&(\GA_V, \Gamma)\ar[r]^{\varphi_{j,V}}\ar[dd]& 	(\GA_{j,V}, \Gamma_j)\ar[dd]^{\alpha}\\
			(\GA,\Gamma)\ar[ru]\ar[rd]&&\\
			&(\GA_U, \Gamma)\ar[r]^{\varphi_{j,U}}&(\GA_{i,U}, \Gamma_i)}$$  
		commutes. Observe that $\alpha$ is unique, if it exists, because $\varphi_{j,V}$ is surjective. Hence $(I, \leq)$ is a partially ordered set. We also observe that the restriction of $\leq$  to $I_U$ coincides with  the partial order of $I_U$ $(U \in \U)$. We claim that this ordering  makes $(I,\leq)$ into a directed poset. 
		
		To see this consider $(j,V), (i,U)\in I$. Let $L\geq V,U$. Then by Lemma \ref{factors through}  there exists $(k,L)\in I_L$ and morphisms of graphs of groups $\underline{\alpha}_i:(\GA_{k,L},\Gamma_k)\longrightarrow  (\GA_{i,U},\Gamma_i)$, $\underline{\alpha}_j:(\GA_{k,L},\Gamma_k)\longrightarrow  (\GA_{j,V},\Gamma_j)$ such that the diagram 
		$$\xymatrix{&(\GA_V, \Gamma)\ar[r]^{\varphi_{j,V}}& 	(\GA_{j,V}, \Gamma_j)\\
			(\GA_L,\Gamma)\ar[ru]\ar[rd]\ar[r]&(\GA_{Lk}, \Gamma_k)\ar[ru]^{\alpha_j}\ar[rd]^{\alpha_i}&\\
			&(\GA_U, \Gamma)\ar[r]^{\varphi_{j,U}}&(\GA_{i,U}, \Gamma_i)}$$
		commutes.
		
		All maps in this diagram are surjective; it follows that $\alpha_i, \alpha_j$ are unique. This shows that $(I,\leq)$ is directed and it follows that $(\GA,\Gamma)=\varprojlim_{U \in \U, i\in  I_U} (\GA_{i,U},\Gamma_i)$. 	
	\end{proof}
	
	\begin{prop}\label{decomposition of graph of groups} Let $(\GA,\G)$ be an injective profinite graph of pro-$\CA$-groups. Then $(\GA,\G)$ decomposes as an inverse limit $(\GA,\G)=\varprojlim_{i \in I}(\GA_i,\Gamma_i)$ of finite  graphs of $\CA$-groups.
		
	\end{prop}	
	
	\begin{proof}  Let $G=\Pi_1^{\CA}(\GA,\Gamma)$ be the fundamental group of $(\GA, \Gamma)$ with respect to  a universal specialization $(\nu,\nu'):(\GA,\Gamma)\longrightarrow \Pi_1^{\CA}(\GA,\Gamma)$. Let $\mathcal{U}$ be the collection of all open normal subgroups of $G$. 
		
		
		Since $\GA_U=\bigcup_{m \in \G} \nu(\GA(m))U/U$ is closed in $G/U$, by \cite[Lemma 5.2.1 combined with page 149]{Ribes} the set $\GA_U=\{(g,m)\in G/U\times\Gamma\mid m\in \Gamma, g\in \nu(\GA(m))U/U\}$ is a sheaf. For simplicity we identify from  now on  the fiber ${\cal G}_U(m)= \nu({\cal
			G}(m))U/U \times \{m\}$ with the subgroup $\nu({\cal G}(m))U/U$ in $G/U$.
		
		Define then a profinite graph of finite groups, $(\GA_U,\G)$  putting $$\partial_j:\GA_U(m)\longrightarrow \GA_U(d_j(m))$$ by $\partial_j(\nu(g)U)=\nu\partial_j(g)U$, for $g\in \GA(m)$, $j=0,1$. To see that it is well-defined let $g,h\in\GA(m)$ such that $\nu(g)U=\nu(h)U$ and let $\chi:\Gamma\longrightarrow \pi_1^{\CA}(\Gamma)$ be the map from Definition \ref{spec}. Then $$\partial_0(\nu(g)U)=\nu\partial_0(g)U=\nu(g)U=\nu(h)U=\nu(\partial_0(h))U=\partial_0(\nu(h)U)$$ and $$\partial_1(\nu(g)U)=\nu\partial_1(g)U=(\nu'\chi(m))^{-1}\nu(g)\nu'\chi(m)U=$$$$=(\nu'\chi(m))^{-1}\nu(h)\nu'\chi(m)U=\nu\partial_1(h)U=\partial_1(\nu(h)U).$$
		
		Put $G_U=\Pi^{\CA}_1(\GA_U,\G)$ and define $\tilde U=\langle U\cap \GA(v)^g \mid v \in V(\Gamma), g\in G\rangle$. It is not difficult to see, using the argument similar to the proof of \cite[Corollary 5.5.9]{Ribes}, that $G_U=G/\tilde U$ and so $G=\varprojlim_{U\triangleleft_o G} G_U$. Thus by Lemma \ref{inverse system} it suffices to show the proposition for $(\GA_U,\G)$. 
		
		
		Let $(\nu_U, \nu'_U):(\GA_U, \Gamma)\longrightarrow G_U$ be the universal specialisation. 
		The projection $G/U\times \Gamma \to G/U$ restricts to a continuous map
		$\rho: {\cal G}_U \to G/U$ that sends ${\cal G}_U (m)$ identically to the
		subgroup $\nu({\cal G}(m))U/U$ of $G/U$.  So $\rho$ induces a unique
		continuous homomorphism $f_U: G_U\to G/U$  such that  $f_U\nu_U= \rho$.
		Hence $f_U\nu_U$ is injective on ${\cal G}_U (m)$,   for each $m\in
		\Gamma$.

		By \cite[Theorem 5.3.4 and Lemma 5.3.1]{Ribes} the sheaf $(\GA_U,\pi_U,\Gamma)$ decomposes as a surjective inverse limit $(\GA_U,\pi_U,\Gamma)=\varprojlim_{i\in I} (\GA_i,\pi_i, \Gamma_i) $ of finite sheaves $(\GA_i,\pi_i, \Gamma_i)$ (the surjectivity follows from the proof); moreover one can achieve  $\Gamma_i$ to be finite graphs by  replacing  in the proof of \cite[Theorem 5.3.4]{Ribes} the set $\mathcal{R}$ of clopen equivalence relations $R$  by the set of clopen equivalence relations such that $\Gamma/R$ is a finite graph.   Let $\underline\varphi_i=(\varphi_i,\varphi'_i):(\GA_U,\Gamma)\longrightarrow (\GA_i, \Gamma_i)$ be the natural projection.    
		
		We need to turn this inverse limit decomposition into a decomposition of graph of groups, i.e to define  the structure of graphs of finite groups on $(\GA_i, \Gamma_i)$ and show that  we have an inverse system of graphs of groups. 
		
		There exists $i_0\in I$ such that the continuous maps $f_U\nu_U,f_U\nu'_U\chi$ factor via $(\GA_{i_0},\Gamma_{i_0})$ (see Lemma \ref{factors through} and \cite[Lemma 1.1.16]{RZ}), i.e. there exist a fiber homomorphism  $\beta_{i_0}:\GA_{i_0}\longrightarrow G/U$ and a continuous map $\beta'_{i_0}:\Gamma_{i_0}\longrightarrow G/U$ such that $f_U\nu_U=\beta_{i_0}\varphi_{i_0}$  and $f_U\nu'_U\chi=\beta'_{i_0}\varphi'_{i_0}$.   
		
		For $i\geq i_0$ define  a fiber homomorphism  $\beta_{i}:\GA_{i}\longrightarrow G/U$ and a continuous map $\beta'_{i}:\Gamma_{i}\longrightarrow G/U$  by 
		$\beta_i=\beta_{i_0}\varphi_{i,i_0},$ , $\beta'_i=\beta'_{i_0}\varphi'_{i,i_0}$, where $$\varphi_{i,i_0}:(\GA_i, \Gamma_i)\longrightarrow (\GA_{i_0}, \Gamma_{i_0}), \varphi'_{i,i_0}:\Gamma_{i}\longrightarrow \Gamma_{i_0}$$ the morphisms of the corresponding projective systems. Then $$f_U\nu_U=\beta_{i}\varphi_i, f_U\nu'_U\chi=\beta'_{i}\varphi'_i \eqno{(*)}$$ for every $i$.
		
		W.l.o.g we  assume from now on that all $i\geq i_0$. In particular, $\beta_{i}$ is injective on $\GA_{i}(m)$ for all $m\in \Gamma_{i}$ and for each $i\in I$.  It follows that the natural morphisms 
		$\varphi_i:\GA_U\longrightarrow \GA_i$ are injective on fibers.
		Thus we have the following commutative diagram:
		$$\xymatrix{\Gamma\ar[d]^{\chi}\ar[r]^{\varphi'_i}&\Gamma_i\ar[dd]^{\beta'_i}\\
			\pi_1^{\CA}(\Gamma)\ar[d]^{\nu'_U}&\\ 
			\Pi_1^{\CA}(\GA_U, \Gamma)\ar[r]^{f_U}& G/U\\
			\GA_U\ar[u]^{\nu_U}\ar[r]^{\varphi_i}& \GA_i\ar[u]^{\beta_i}}$$
		
		\medskip
		
		Define $\partial_j\varphi_i(g)=\varphi_i\partial_j(g)$, $g\in \GA_i(m)$ on $(\GA_i,\Gamma_i)$, $(j=0,1)$. 
		We  show that $\partial_j$ is well-defined on $ (\GA_i, \Gamma_i)$, $(j=0,1)$.  Choose $g\in \GA_U(e)$, $h\in \GA(e')$ such that $\varphi'_i(e)=\varphi'_i(e')$ and $\varphi_i(g)=\varphi_i(h)$.
		Since $f_U\nu_U$ is injective on $\GA_U(m)$ and  $\beta_i\varphi_i(x)=f_U\nu_U(x)$ for every $x\in \GA_U$ by $(*)$, it suffices to show that 	$ f_U	\nu_U\partial_j(g)=f_U\nu_U\partial_j(h)$ for $j=0,1$. 
		Since $\nu_U(x)=\nu_U\partial_0(x)$ for any $x\in \GA_U$ we have $$ 	f_U\nu_U\partial_0(g)=f_U\nu_U(g)=\beta_{i}\varphi_i(g)=\beta_{i}\varphi_i(h)=f_U\nu_U(h)=f_U\nu_U\partial_0(h);\eqno{(**)}$$ so  $\partial_0$ is well-defined.  Now  we need to show that  	$ f_U	\nu_U\partial_1(g)=f_U\nu_U\partial_1(h)$. But $$ f_U	\nu_U\partial_1(g)=f_U( (\nu'_U\chi(e))^{-1}\nu_U(g)\nu'_U\chi(e))  =f_U (\nu'_U\chi(e)^{-1})f_U(\nu_U(g))f_U\nu'_U\chi(e)\overset{(*)}{=}$$ 
		$$\overset{(*)}{=} (\beta'_i\varphi'_i(e))^{-1}f_U\nu_U(g)\beta'_i\varphi'_i(e)\overset{(**)}{=} (\beta'_i\varphi'_i(e))^{-1}f_U\nu_U(h)\beta'_i\varphi'_i(e)=$$$$=(\beta'_i\varphi'_i(e'))^{-1}f_U\nu_U(h)\beta'_i\varphi'_i(e')\overset{(*)}{=}  f_U (\nu'_U\chi(e'))^{-1}\nu_U(h)\nu'_U\chi(e'))=   f_U\nu_U\partial_1(h)$$  as desired. 
		
		It remains to observe that $\partial_i$ are injective on $(\GA_i,\Gamma_i)$ since $f_U\nu_U$ is injective on $\GA_U(m)$ for all $m\in \Gamma$.
		
		Thus $(\GA_U,\Gamma)=\varprojlim_i (\GA_i, \Gamma_i)$ is a decomposition as an inverse limit of finite graphs of finite $\CA$-groups as required.
	\end{proof}
	
	\begin{obs}\label{converse} The hypothesis of injectivity of $(\GA,\Gamma)$ 
		in Proposition \ref{decomposition of graph of groups} is essential. Indeed, if $\CA$ consists of all finite groups,  then $\Pi_1^{abs}(\GA_i,\G_i)$ is virtually free and so $\Pi_1(\GA_i,\G_i)=\widehat{\Pi_1^{abs}(\GA_i,\G_i)}$; thus $(\GA_i,\Gamma_i)$ is injective for each $i\in I$.  It follows that $(\GA,\Gamma)=\varprojlim_i (\GA_i, \Gamma_i)$ is injective.
	\end{obs}
	
	Combining 	Proposition \ref{decomposition of graph of groups} and Remark \ref{converse} we can state the following
	
	\begin{cor} Let $(\GA,\G)$ be a  profinite graph of profinite groups. Then $(\GA,\G)$ decomposes as an inverse limit $(\GA,\G)=\varprojlim_{i \in I}(\GA_i,\Gamma_i)$ of finite  graphs of finite groups if and only if $(\GA,\G)$ is injective.
		
	\end{cor}
	
	We shall recall now the definitions the standard tree and standard pro-$\CA$ tree on which the fundamental group naturally acts following \cite{Serre} and \cite{Ribes} respectively; we shall need them in the following sections. 	
	
	\begin{df}[Standard tree] Let $G=\pi_1(\GA,\G)$ be the fundamental group of the graph of abstract groups $(\GA,\G)$. There is an abstract standard graph $S^{abs}=S^{abs}(\GA,\G)$ which is in fact a tree (cf. \cite{Serre}, Sec. I.5.3). Let $T$ be a maximal subtree of $\G$. We define $S^{abs}$ by
		\begin{center}
			$V(S^{abs})=\bigcup_{v \in V(\G)} G/\GA(v)$ and $E(S^{abs})=\bigcup_{e \in E(\G)} G/\GA(e)$
		\end{center}
		and its incidence maps
		\begin{center}
			$d_0(g\GA(e))=g\GA(d_0(e))$, $d_1(g\GA(e))=gt_e\GA(d_1(e))$
		\end{center}
		$(g \in G, e \in E(\G))$ and $t_e=1, \forall e \in E(T)$.
	\end{df}	
	
	Next one defines the $\CA$-standard graph associated with a graph of pro-$\CA$ groups.

	\begin{df}[$\CA$-standard graph of a graph of pro-$\CA$ groups] \label{cstd}
		Let $(\GA,\pi,\G)$ be a graph of pro-$\CA$ groups over a connected profinite graph $\G$, $j:\G \ra \widetilde{\G}$ be a continuous 0-section of the universal Galois $\CA$-covering $\zeta:\widetilde{\G} \ra \G$ of $\G$, and let $J=j(\G)$ be the corresponding 0-transversal. Let $(\gamma,\gamma')$ be a $J$-specialisation of $\GG$ in the fundamental pro-$\CA$ group of $(\GA,\pi,\G)$, $\Pi=\Pi_1^{\CA}(\GA,\pi,\G)$. Then we can define a profinite graph $S=S^{\CA}(\GA,\G)= S^{\CA}(\GA,\G,\Pi)$ which is canonically associated to the graph of groups $\GG$ and $\FG$. 
		
		For $m \in \G$, define $\Pi(m)=\gamma(\GA(m))$. As a topological space, $S^{\CA}(\GA,\G)$ is defined to be the quotient space of $\G \times \Pi$ modulo the equivalence relation $\sim$ given by 
		\begin{center}
			$(m,h) \sim (m',h')$ if $m=m',h^{-1}h' \in \Pi(m)$ $(m,m' \in \G, h,h' \in \Pi).$
		\end{center}
		So, as a set, $S^{\CA}(\GA,\G)$ is the disjoint union $$S^{\CA}(\GA,\G)=\bigcup_{m \in \G} \Pi/\Pi(m).$$ Denote by $\alpha:\G \times \Pi \ra S^{\CA}(\GA,\G)$ the quotient map. The projection $p':\G \times \Pi \ra \G$ induces a continuous epimorphism $p:S^{\CA}(\GA,\G) \ra \G,$ such that $p^{-1}(m)=\Pi/\Pi(m)$ and $p'=p \alpha$.
		
		To make $S^{\CA}(\GA,\G)$ into a profinite graph we define the subspace of vertices of $\SG$ by $V(\SG)=p^{-1}(V(\G))$ and the incidence maps by
		\begin{equation*} \label{4.3}
			d_0(hH(m))=hH(d_0(m))
		\end{equation*}
		\begin{equation*} \label{6.8}
			d_1(hH(m))=h(\gamma'\chi(m))H(d_1(m)),
		\end{equation*}
		$(h \in H, m \in \G)$. The definition of $S^{\CA}(\GA,\G)$ is independent, up to isomorphism, on the choice of the 0-section $j$ (cf. \cite{Ribes}, Theorem $6.3.3$)
	\end{df}
	
	There is a natural continuous action of $\Pi=\FG$ on the graph $S^{\CA}(\GA,\G)$ given by \[g(g'\Pi(m))=gg'\Pi(m),\] $(g,g' \in \Pi, m \in \G)$. 
	
	The $\CA$-standard graph  is in fact a $\CA$-tree when $\CA$ is a pseudovariety of finite groups which is extension closed (cf. \cite{Ribes}, Corollary 6.3.6). 
	
	\section{Definition of the fundamental pro-$\CA$ group at based point}

	The  goal of this section is to define
	the pro-$\CA$ fundamental group
	of an injective graph of pro-$\CA$ groups via an inverse limit of completions of the
	abstract fundamental group of finite  graphs of finite $\CA$-groups with
	respect to a point in the original graph. This is a new concept in the profinite version of Bass-Serre theory and has the advantage of behaving better with repect to morphisms of finite graph of finite groups and even profinite graphs of pro-$\CA$ groups.  One can not use classical definition of the fundamental group of graph of groups with respect to a maximal subtree, because  the image of a maximal subtree under epimorphism of graphs is not a maximal subtree	in general.
	
	\begin{df}[The group $F(\GA,\G)$ {\cite[Sect.\, I.5.1]{Serre}}]
		The path group $F(\GA,\G)$ is defined by $F(\GA,\G)=W_1/N$, where $W_1=\left(\mathop{\Asterisk}_{v \in V(\G)} \GA(v)\right) \Asterisk F(E(\G))$, where $F(E(\G))$ denotes the free group with basis $E(\G)$ and $N$ is a normal subgroup of $W_1$ generated by the set $\left\{ \partial_0(x)^{-1}e\partial_1(x)e^{-1} \mid x \in \GA(e), e \in E(\G) \right\}$. 
	\end{df}
	
	\begin{df}[Words of $F(\GA,\G)$ {\cite[Sec.\,I.5.1, Definition 9]{Serre}}]
		Let $c=v_0,e_0, \cdots, e_n,v_n$, be a path in $\G$ with length $n=l(c)$ such that $v_j \in V(\G), e_j \in E(\G)$, $j=0, \cdots, n$. A word of type $c$ in $F(\GA,\G)$ is a pair $(c,\mu)$ where $\mu=(g_0, \cdots, g_n)$ is a sequence of elements $g_j \in \GA(v_j)$. The element $|c,\mu|=g_0,e_0,g_1,e_1, \cdots,e_n,g_n$ of $F(\GA,\G)$ is said to be associated with the word $(c,\mu)$.   
	\end{df}	
	
	\begin{df}[The fundamental group of $(\GA,\G)$ {\cite[Sect. I.5.1, Definition 9(a)]{Serre}}]
		Let $v$ be a vertex of $\G$. We define $\pi_1(\GA,\G,v)$ as the set of elements of $F(\GA,\G)$ of the form $|c,\mu|$, where $c$ is a path whose extremities both equal $v$. One sees immediately that $\pi_1(\GA,\G,v)$ is a subgroup of $F(\GA,\G)$, called the fundamental group of $(\GA,\G)$ at $v$. In particular, if $\GA$ consists of trivial groups only then  $\pi_1(\GA,\G,v)$ becomes a usual fundamental group of the graph $\Gamma$ and denoted by  $\pi_1(\Gamma, v)$. It can be viewed of course as a subgroup that consists of set of elements of $F(\GA,\G)$ of the form $|c,\mu|=g_0,e_0,g_1,e_1, \cdots,e_n,g_n$, where $c$ is a path whose extremities both equal $v$ and $g_0=1=g_1=\ldots =g_n$. This way $G=\pi_1(\GA,\G,v)$ is a semidirect product $\pi_1(\GA,\G,v)=\langle \GA(v)\mid v\in V(\Gamma)\rangle^G \rtimes \pi_1(\Gamma, v)$.
	\end{df}

	\begin{df} If $\Gamma$ is a connected finite graph,	its pro-$\CA$ fundamental group $\pi_1^{\CA}(\Gamma,v)$ can be defined as the pro-$\CA$ completion $\pi_1(\Gamma, v)_{\widehat \CA}$ of $\pi_1(\Gamma, v)$. If $\Gamma$ is a connected profinite graph and $\Gamma=\li \Gamma_i$ its decomposition as an inverse limit of  finite graphs $\Gamma_i$, then $\pi_1^{\CA}(\Gamma, v)$ can be defined as the inverse limit  $\pi_1^{\CA}(\Gamma)=\li \pi_1^{\CA}(\Gamma_i, v_i )$, where $v_i$ is the image of $v$ in $\Gamma_i$ (see \cite[Proposition 3.3.2 (b)]{Ribes}). \end{df}
	
	Next we shall use a similar approach to define the fundamental pro-$\CA$ group  $\Pi_1^{\CA}(\GA,\Gamma, v)$ at based point $v$.
	
	\begin{prop}[{\cite[Proposition $2.4$]{Bass}}] \label{ind} A morphism of graphs of groups $$\underline{\alpha}=(\alpha,\alpha'):(\GA,\G) \ra (\GA',\G')$$ induces a morphism of fundamental groups $$\beta:\pi_1(\GA,\G, v) \ra \pi_1(\GA',\G',v')$$ defined by $\beta(|c,\mu|)=|\alpha'(c),\alpha(\mu)|$, where $\alpha(\mu)=(\alpha(g_0),\alpha(g_1), \cdots, \alpha(g_n))$, $g_i \in \GA(v_i)$ and $\beta(v)=\alpha'(v)=v'$. 
	\end{prop}
	
	\begin{prop} \label{fg}
		An inverse limit $(\GA,\G)=\li_{i \in I} (\GA_i,\G_i)$ of finite abstract graphs of finite groups (that are in $\CA$) induces an inverse limit $\li_{i \in I} (\pi_1(\GA_i,\G_i,v_{i}))_{\widehat{\CA}}$ of the pro-$\CA$ completions of fundamental abstract groups $ \pi_1(\GA_i,\G_i,v_{i})$. 
	\end{prop}
	\begin{proof}
		Let $\{(\GA_i,\G_i),\underline{\alpha}_{ij},I\}$ be the corresponding inverse system of finite abstract graphs of groups. By Proposition \ref{ind}, the morphism $\underline{\alpha}_{ij}: (\GA_i,\G_i) \ra (\GA_j,\G_j)$ induces a morphism of fundamental groups $\beta_{ij}:\pi_1(\GA_i,\G_i, v_{i}) \ra \pi_1(\GA_j,\G_j,v_{j})$ defined by $\beta(|c,\mu|)=|\alpha'(c),\alpha(\mu)|$, where $\alpha(\mu)=(\alpha(g_0),\alpha(g_1), \cdots, \alpha(g_n))$, $g_j \in G(v_j)$ and $\beta(v_{i})=\alpha'(v_{i})=v'_{i}$. 
		
		Therefore, we can construct an inverse system $\{(\pi_1(\GA_i,\G_i, v_{i}))_{\widehat{\CA}},\beta_{ij},I\}$.  
	\end{proof}
	
	We are ready to define the pro-$\CA$ fundamental group $\Pi_1^{\CA}(\GA,\G,v)$ of a profinite graph of pro-$\CA$ groups with base vertex $v$ of $\Gamma$. Note that it works for an injective profinite graph of pro-$\CA$-groups only (see Remark \ref{converse}. 
	
	\begin{df} \label{fundamental group at point} Let $(\GA,\Gamma)$ be an injective profinite graph of pro-$\CA$ groups and $(\GA,\G)=\li_{i \in I} (\GA_i,\G_i)$ be the decomposition as the inverse limit of  finite  graphs of finite  $\CA$-groups (See Proposition \ref{decomposition of graph of groups}). Let $v$ be a vertex of $\Gamma$. The group $\li_{i \in I} (\pi_1(\GA_i,\G_i,v_{i}))_{\widehat{\CA}}$ from Proposition \ref{fg}  will be called the pro-$\CA$ fundamental group of the graph of pro-$\CA$ groups $(\GA, \Gamma)$ at point $v$ and denoted by $\Pi_1^{\CA}(\GA,\G,v)$.\end{df}
	
	
	
	The next proposition shows that	the definition of $\Pi_1^{\CA}(\GA,\G,v)$ does not depend on  decomposition of $(\GA,\G)$  as an inverse limit of finite graphs of finite groups.
	
	\begin{prop}\label{decomposition independence} 
		Let $(\GA,\G)$ be an injective profinite graph of pro-$\CA$-groups and $(\GA,\G)=\li_{i \in I} (\GA_i,\G_i)$ be a decomposition as an inverse limit of  finite  graphs of finite  $\CA$-groups. Then $$\Pi_1^{\CA}(\GA,\G,v)=\li_{i \in I} \Pi_1^{\CA}(\GA_i,\G_i,v_{i}).$$
	\end{prop}
	
	\begin{proof} Let $\underline\eta:(\GA, \Gamma)\longrightarrow (\HA, \Delta)$ be an epimorphism to a finite  graph of $\CA$-groups $(\HA, \Delta)$. Choose $v\in V(\Gamma)$ and let $v_0$ be its image in $\Delta$. It suffices to show that the natural epimorphism $\Pi_1^{\CA} (\GA, \Gamma, v)\longrightarrow \Pi_1^{\CA}(\HA, \Delta, v_0)$ factors through some $\Pi_1^{\CA}(\GA_i,\G_i,v_{i})$. 
		
		Since $(\HA, \Delta)$ is finite, by Lemma \ref{factors through} $\underline\eta$ factors through some $\underline\eta_i:(\GA_i,\G_i)\longrightarrow (\HA, \Delta)$, i.e. $\underline\eta=\underline\eta_i\underline\pi_i$, where $\underline\pi_i:(\GA, \Gamma)\longrightarrow (\GA_i,\G_i)$ is the natural projection. Let $v_i$ be the image of $v$ in $\Gamma_i$. By Proposition \ref{ind} $\underline\eta_i$ induces the natural homomorphism $\pi_1(\GA_i,\G_i, v_i)\longrightarrow \pi_1(\HA, \Delta, v_0)$ that in turn induces  the homomorphism of pro-$\CA$ completions $$\Pi_1^{\CA}(\GA_i,\G_i,v_i)=(\pi_1(\GA_i,\G_i, v_i))_{\widehat\CA}\longrightarrow (\pi_1(\HA, \Delta, v_0))_{\widehat\CA}=\Pi_1^{\CA}(\HA, \Delta, v_0).$$ Thus we have a diagram 
		$$\xymatrix{\Pi_1^{\CA} (\GA, \Gamma, v)\ar[dr]\ar[d]&\\
			\Pi_1^{\CA}(\GA_i,\G_i,v_i)\ar[r] &\Pi_1^{\CA}(\HA, \Delta, v_0)}$$  that commutes on vertex groups and underlying graphs. Hence this diagram commutes and the proof is finished. 
		
	\end{proof}

	We shall show now that for injective $(\GA,\Gamma)$ this definition and Definition \ref{fgg} are equivalent.
	
	\begin{theo}\label{equivalence} Let $(\GA,\Gamma)$ be an injective profinite graph of pro-$\CA$ groups. Then $\Pi_1^{\CA}(\GA, \G, v)\cong \Pi_1^{\CA}(\GA, \Gamma)$.	
	\end{theo}
	
	\begin{proof} Let $\widetilde{\G}$ be the universal $\CA$-covering of the profinite graph $\G$ and $\zeta: \widetilde{\G} \ra \G$ be the covering map. 
		
		Consider  a graph of groups $(\widetilde{\GA},\widetilde{\G})$ from Definition \ref{universal graph of groups}. Choose a base point $v\in \Gamma$. The action of $\pi_1^{\CA}(\Gamma,v)$ on $\widetilde{\G}$ induces the action on  $(\widetilde{\GA},\widetilde{\G})$ that in turn induces the action of $\pi_1^{\CA}(\Gamma,v)$ on $\Pi_1^{\CA}(\widetilde{\GA},\widetilde{\G})$, i.e. we have a semidirect product $\Pi_1^{\CA}(\widetilde{\GA},\widetilde{\G})\rtimes \pi_1(\Gamma, v)$. 
		
		Choose  $\tilde v\in \widetilde \Gamma$ such that $\zeta(\tilde v)=v$. Let $(\GA,\Gamma)=\li_i (\GA_i,\G_i)$ be a decomposition of $\Gamma$ as the inverse limit of finite quotient graphs of groups $(\GA_i,\Gamma_i)$ (cf. Proposition \ref{decomposition of graph of groups}). By Proposition 3.3.2 \cite{Ribes} the decomposition $\Gamma=\li_i \G_i$ induces an inverse system of  pairs of compatible morphisms  $\zeta_i: \widetilde \Gamma_i\rightarrow  \Gamma_i$, $f_{ij}:\pi_1(\Gamma_i)\rightarrow \pi_1(\Gamma_j)$ such that  $\pi_1^{\CA}(\Gamma)=\li_i \pi_1^{\CA}(\Gamma_i, v_i)$, $\zeta=\li_i \zeta_i$, where $v_i$ is the image of $v$ in $\Gamma_i$. This defines a decomposition as the inverse limit $(\widetilde \GA,\widetilde \Gamma)=\li_i (\widetilde\GA_i,\widetilde \G_i)$ (as $(\widetilde\GA_i,\widetilde \G_i)$ is defined as a pull-back of $\GA_i\longrightarrow \Gamma_i$ and $\zeta_i:\widetilde \Gamma_i\longrightarrow \Gamma$, see Definition \ref{universal graph of groups}).
		
		Let $\tilde v_i$ be the image of $\tilde v$ in $\widetilde\Gamma_i$. Let $c=v_i=v_0,e_0, \cdots, e_n,v_n=v_i$, be a circuit of length $n$ such that $v_j \in V(\G), e_j \in E(\G)$, $j=0, \cdots, n$ that we regard also as an element of $\pi_1(\Gamma_i, v_i)$.  Let $\widetilde\Gamma_i^{abs}$ be the  connected component of $\widetilde\Gamma_i$ (regarded as an abstract graph) containing $\tilde v_i$. Then $\widetilde\Gamma_i^{abs}$ is the usual universal cover of $\Gamma_i$ (see \cite[Proposition 8.2.4]{Ribes}).  It follows that the  circuit $c$ lifts to the unique path $\tilde c=\tilde v_i=\tilde v_0,\tilde e_0, \ldots, \tilde e_n,\tilde v_n=c\tilde v_i$ from $\tilde v_i$ to $c\tilde v_i$ in $\widetilde\Gamma_i^{abs}$. 
		
		Denote by $(\widetilde{\GA}_i,\widetilde\G_i^{abs})$ the graph of groups obtained by the restriction of $(\widetilde{\GA}_i,\widetilde{\G_i})$ to $\widetilde \G_i^{abs}$ and let $\pi_1(\widetilde{\GA}_i,\widetilde\G_i^{abs}, \tilde v_i)$ be its fundamental group. Then $\pi_1(\Gamma_i, v_i)$ acts naturally on $(\widetilde{\GA}_i,\widetilde\G_i^{abs})$ that induces the action  of $\pi_1(\Gamma_i, v_i)$ on  $\pi_1(\widetilde{\GA}_i,\widetilde\G_i^{abs}, \tilde v_i)$ so that we can consider the semidirect product $\pi_1(\widetilde{\GA}_i,\widetilde\G_i^{abs}, \tilde v_i)\rtimes \pi_1(\Gamma_i, v_i)$.
		
		An element $|c,\mu|=g_0,e_0,g_1,e_1, \cdots,e_n,g_n$  in $ \pi_1(\GA_i,\G_i,v_i)$, where $\mu=(g_0, \cdots, g_n)$ is a sequence of elements $g_j \in \GA(v_j)$, lifts to a unique element $(\tilde{\mu}, c)\in \pi_1(\widetilde{\GA}_i,\widetilde{\G}^{abs}_i, \tilde v_i)\rtimes \pi_1(\Gamma_i, v_i)$, where $\tilde\mu=\tilde g_0 \tilde g_1 \cdots \tilde g_n\in \pi_1(\widetilde{\GA},\widetilde{\G}, \tilde v_i)$ with $\tilde g_j=\tilde\zeta_{|\widetilde \GA(\tilde v_j)}^{-1}(g_j)$.
		
		Hence we can define a map $\psi_i: \pi_1(\GA_i,\G_i,v_i) \ra \pi_1(\widetilde{\GA_i},\widetilde{\G_i}^{abs}, \tilde v_i) \rtimes \pi_1(\G_i, v_i)$ by $\psi_i(|c,\mu|)=(\tilde{\mu},c)$. To see that $\psi_i$ is a homomorphism, let $|c,\mu|=|c_1,\mu_1||c_2,\mu_2|$ and let $\tilde c_1$, $\tilde c_2$ be liftings of $c_1$, $c_2$ respectively. Then $\psi_i(|c,\mu|)=(\tilde\mu, c)=(\tilde\mu_1(c_1\tilde \mu_2c_1^{-1}), c_1c_2)$ as needed.    
		
		Thus $\psi_i$ is a homomorphism which is  is clearly bijective, i.e. $\psi_i$ is an isomorphism. Moreover, for $i> j$ the commutative diagram
		
		\begin{equation*}
			\begin{tikzcd}
				\widetilde\G_j\ar[d,"\zeta_j",swap] \ar[r]&\widetilde\G_i\ar[d, "\zeta_i"]\\ 
				\G_j\ar[r]& \G_i
			\end{tikzcd}
		\end{equation*}
		
		induces the commutative diagram 
		
		\begin{equation*}
			\begin{tikzcd}
				\pi_1(\GA_j, \G_j, v_j)\ar[d,"\psi_j",swap] \ar[r]&\pi_1(\GA_i, \G_i, v_i)\ar[d,"\psi_i"]\\ \pi_1(\widetilde{\GA}_j,\widetilde{\G_j}^{abs}, v_j) \rtimes \pi_1(\G_j, v_j)\ar[r]& \pi_1(\widetilde\GA_i, \widetilde\G_i^{abs}, v_i)\rtimes \pi_1(\G, v_i)
			\end{tikzcd}
		\end{equation*}
		
		\medskip
		Now recall from Definition \ref{fundamental group at point} that $\Pi_1^{\CA}(\GA_i,\G_i,v_i)=\pi_1(\GA_i,\G_i,v_i)_{\widehat \CA}$. It is straightforward to check that the pro-$\CA$ completion of $\pi_1(\widetilde{\GA}_i,\widetilde{\G_i}^{abs}, \tilde v_i) \rtimes \pi_1(\G_i, v_i)$  gives $\Pi_1^{\CA}(\widetilde{\GA_i},\widetilde{\G_i}) \rtimes \pi_1^{\CA}(\G_i, v_i)$ as defined above. This shows that $\psi$ induces the isomorphism  $$\hat\psi: \Pi_1^{\CA}(\GA_i,\G_i,v_i) \ra \Pi_1^{\CA}(\widetilde{\GA_i},\widetilde{\G_i}) \rtimes \pi_1^{\CA}(\G, v_i)$$ and we have the following commutative diagram 
		
		\begin{equation*}
			\begin{tikzcd}
				\Pi_1^{\CA}(\GA_j, \G_j, v_j)\ar[d,"\hat\psi_j",swap] \ar[r]&\Pi_1^{\CA}(\GA_i,\G_i,v_i)\ar[d, "\hat\psi_i"]\\ \Pi_1^{\CA}(\widetilde{\GA}_j,\widetilde{\G_j}, \tilde v_j) \rtimes \pi_1^{\CA}(\G_j, v_j)\ar[r]& \Pi_1^{\CA}(\widetilde\GA_i, \widetilde\G_i, \tilde v_i)\rtimes \pi_1^{\CA}(\G, v_i)
			\end{tikzcd}
		\end{equation*}
		
		Then $\psi=\varprojlim_{i} \hat\psi_i$ is an isomorphism $\psi: \Pi_1^{\CA}(\GA, \G, v)\longrightarrow  \Pi_1^{\CA}(\widetilde{\GA},\widetilde{\G}, \tilde v) \rtimes \pi_1^{\CA}(\G, v)$. 
		
		Note that $\Pi_1^{\CA}(\GA, \Gamma) \cong \Pi_1^{\CA}(\widetilde \GA, \widetilde\Gamma)\rtimes \pi_1^{\CA}(\Gamma, v)$ by \cite[Proposition 6.5.1]{Ribes}). Therefore  $\Pi_1^{\CA}(\GA, \G, v)\cong \Pi_1^{\CA}(\GA, \Gamma)$ as desired.	
	\end{proof}

	\section{The pro-$\CA$ completion of the fundamental group of an infinite abstract graph of groups} 
	
	The profinite version of the Bass-Serre theory can be used very effectively in the study of certain abstract groups. One sees such a group $G$  as the fundamental group of a graph of groups and then the profinite completion $\widehat G$ is the fundamental group $\widehat G=\Pi_1(\widehat \GA,\Gamma)$ of the same graph of the profinite completions of edge and vertex groups.  Through this view then it is possible to apply geometric techniques to obtain algebraic results. Until now, one could use  this approach only for finitely generated groups, assuming that the graph $\G$ in the graph of groups $(\GA, \G)$ is finite. The main reason for this assumption is that then $\G$ is automatically a profinite graph and so $(\widehat{\GA},\G)$ is automatically a profinite graph of pro-$\CA$-groups, where each $\widehat{\GA}(m)$ is the pro-$\CA$ completion of $\GA$. One also has a natural way of relating $\pi_1(\GA,\G)$ with $\Pi=\Pi_1(\widehat{\GA},\G)$ and $S^{abs}=S^{abs}(\GA,\G)$ with $S=S(\widehat{\GA},\G)$. With the aim to apply this technique to abstract groups which are not necessarily finitely generated, we  construct  a profinite graph of pro-$\CA$ groups $(\overline{\GA},\overline{\G})$ where an infinite abstract graph $\G$ in densely embedded in $\overline{\G}$ and the above properties are preserved, the content of Theorem \ref{C}. We assume  that $(\GA,\Gamma)$ is reduced, i.e. whenever $e$ is an edge of $\G$ which is not a loop, then $\partial_i(\GA(e))$ is a proper subgroup of $\GA(d_i(e))$ $(i=0,1)$. This of course does not affect  generality since for any graph of groups we can collapse fictitious edges successfully (i.e. the edges not loops such that $\partial_i(\GA(e))=\GA(d_i(e))$ for  $i=0$ or,$1$) to arrive to the reduced graph of groups.
	
	We finish this section with the proof of Theorem \ref{B} which answers Open Question $6.7.1$ of \cite{Ribes}.
	
	Parts $(a)$, $(b)$ and $(c)$ of Theorem \ref{C} is the subject of  the following 
	
	\begin{theo} \label{thm}
		Let $(\GA,\G)$ be a reduced graph of  groups and $G=\pi_1(\GA,\G)$  its fundamental group. Assume that $G$ is residually $\CA$ and denote by $\overline{\GA}(m)$ the closure of $\GA(m)$ in $G_{\widehat\CA}$. Then 
		\begin{enumerate}[(a)]
			\item There exists an injective profinite graph of pro-$\CA$ groups $(\overline{\GA},\overline{\G})$ such that $\G$ is densely embedded in $\overline{\G}$;
			\item for each $m\in\Gamma$ its vertex group is $\overline{\GA}(m)$;
			\item The fundamental pro-$\CA$ group $\Pi=\Pi_1^{\CA}(\overline{\GA},\overline{\G})$ of $(\overline{\GA},\overline{\G})$ is the pro-$\CA$ completion of $G$ and  so $(\overline{\GA},\overline{\G})$ is injective.
			
			\item  The graph of groups $(\overline{\GA},\overline{\G})$ decomposes as a surjective inverse limit $(\overline{\GA},\overline{\G})=\li (\GA_U,\G_U)$ of finite graphs of finite $\CA$-groups and  
			$
			\Pi_1^{\CA}(\overline{\GA},\overline{\G}, v)=\li (\pi_1(\GA_U,\G_U, v_U))_{\widehat{\CA}}$.
		\end{enumerate}
	\end{theo}
	
	\begin{proof}
		Let $\mathcal{U}$ be the collection of all open normal subgroups of $G$ in the pro-$\CA$ topology of $G$. For $m \in \G$, $U \in \mathcal{U}$, define $\GA_U(m)=\GA(m)U/U$. As $\GA_U(m) \leq G/U$, one concludes that each $\GA_U(m) \in \CA$. 
		
		Define the profinite space $\GA_U=\bigcup_{m \in \G} \GA_U(m)$  and $R_U$ to be the following equivalence relation in $\G$: given $v,w \in V(\G)$, $v \sim_{R_U} w$ if $\GA_U(v)=\GA_U(w)$ and given $e,e' \in E(\G)$, $e \sim_{R_U} e'$ if $d_0(e) \sim_{R_U} d_0(e')$ and $d_1(e) \sim_{R_U} d_1(e')$. Hence, the quotient graph $\G_U$ defined by $\G_U=\G/R_U$ is finite.
		
		Define then a finite graph of finite groups, $(\GA_U,\G_U)$ by putting $\GA_U(\bar m)=\GA_U(m)$ where $\bar m$ is the equivalence class of $m$ and defining $\partial_i(gU/U)=\partial_i(g)U/U$. To see that $\partial_i$ are well-defined, one follows precisely the argument of the third paragraph of the proof of Proposition \ref{decomposition of graph of groups} ignoring the topology.
		
		Put $G_U=\pi_1(\GA_U,\G_U)$ and, for an open subgroup  $Y \leq_o U$, $\underline{\alpha}_{YU}=(\alpha_{YU},\alpha'_{YU}):(\GA_Y,\G_Y) \ra (\GA_U,\G_U)$ to be the epimorphism of graphs of groups defined by $\alpha(gY)=gU$ and $\alpha'(mR_Y)=mR_U$, for $m \in \G$ and $g \in \GA(m)$. These are the canonical quotient maps, so we have the following commutative diagram:
		
		\begin{equation*}
			\begin{tikzpicture}
				\node (0) at (0,0) {$(\GA_Y,\G_Y)$};
				\node (1) at (3,1) {$(\GA_U,\G_U)$};
				\node (2) at (3,-1) {$(\GA_Z,\G_Z)$};
				\node (3) at (6,0) {$(\GA_W,\G_W)$};
				
				\draw[->] (0) edge node[left,sloped, above, pos=0.5] {$\underline{\alpha}_{YU}$} (1);
				\draw[->] (0) edge node[left,sloped, below, pos=0.5] {$\underline{\alpha}_{YZ}$} (2);
				\draw[->] (1) edge node[left,sloped, above, pos=0.5] {$\underline{\alpha}_{UW}$} (3);
				\draw[->] (2) edge node[left,sloped, below, pos=0.5] {$\underline{\alpha}_{ZW}$} (3);
			\end{tikzpicture}
		\end{equation*}
		Indeed,
		$(\alpha_{UW} \circ \alpha_{YU})(gY)=\alpha_{UW}(\alpha_{YU}(gY))=\alpha_{UW}(gU)=gW$ and $(\alpha_{ZW} \circ \alpha_{YZ})(gY)=\alpha_{ZW}(\alpha_{YZ}(gY))=\alpha_{ZW}(gZ)=gW$. For the graph maps, $(\alpha_{UW}' \circ \alpha_{YZ}')(mR_Y)=\alpha_{UW}'(\alpha_{YU}'(mR_Y))=\alpha_{UW}'(mR_U)=mR_W$ and $(\alpha_{ZW}' \circ \alpha_{YZ}')(mR_Y)=\alpha_{ZW}'(\alpha_{YZ}'(mR_Y))=\alpha_{ZW}'(mR_Z)=mR_W$.  
		Hence
		\begin{equation} \label{com}
			\underline{\alpha}_{UW} \circ \underline{\alpha}_{YU}=\underline{\alpha}_{ZW} \circ \underline{\alpha}_{YZ}
		\end{equation}
		for every $m \in \G$ and $g \in \GA(m)$.
		
		Choose a vertex $v$ in $\Gamma$ and denote by $v_U$ its image in $\Gamma_U$. By  Proposition \ref{ind}  this diagram  induces the diagram of the fundamental groups  
		
		\begin{equation*}
			\begin{tikzpicture}
				\node (0) at (0,0) {$\pi_1(\GA_Y,\G_Y, v_Y)$};
				\node (1) at (3,1) {$\pi_1(\GA_U,\G_U,v_U)$};
				\node (2) at (3,-1) {$\pi_1(\GA_Z,\G_Z, v_Z)$};
				\node (3) at (6,0) {$\pi_1(\GA_W,\G_W, v_W)$};
				
				\draw[->] (0) edge node[left,sloped, above, pos=0.5] {$\beta_{YU}$} (1);
				\draw[->] (0) edge node[left,sloped, below, pos=0.5] {$\beta_{YZ}$} (2);
				\draw[->] (1) edge node[left,sloped, above, pos=0.5] {$\beta_{UW}$} (3);
				\draw[->] (2) edge node[left,sloped, below, pos=0.5] {$\beta_{ZW}$} (3);
			\end{tikzpicture}
		\end{equation*}
		Given $|c,\mu| \in \pi_1(\GA_Y,\G_Y)$,  we have $\beta_{UW}(\beta_{YU}(|c,\mu|))=\beta_{UW}(|\alpha_{YU}'(c),\alpha_{YU}(\mu)|)=|(\alpha_{UW}' \circ \alpha_{YU}')(c),(\alpha_{UW} \circ \alpha_{YU})(\mu)|$ and $\beta_{ZW}(\beta_{YZ}(|c,\mu|))=\beta_{ZW}(|\alpha_{YZ}'(c),\alpha_{YZ}(\mu)|)=|(\alpha_{YZ}' \circ \alpha_{YZ}')(c),(\alpha_{ZW} \circ \alpha_{YZ})(\mu)|$, which implies, by Equation (\ref{com}), that 
		\begin{equation} 
			\beta_{UW} \circ \beta_{YU}=\beta_{ZW} \circ \beta_{YZ}
		\end{equation}
		for every $|c,\mu| \in \pi_1(\GA_Y,\G_Y)$. Therefore the diagram commutes.  
		
		Thus we have an inverse system of graphs of groups $\{(\GA_{U},\G_{U}),\underline{\alpha}_{YU}\}$. Put $(\overline{\GA},\overline{\G})=\li (\GA_U,\G_U)$. By Proposition \ref{decomposition independence}, we have 
		\begin{equation} \label{comp}
			\Pi_1^{\CA}(\overline{\GA},\overline{\G}, v)=\li (\pi_1(\GA_U,\G_U, v_U))_{\widehat{\CA}}.
		\end{equation}
		
		Note that  by equation \eqref{comp} $\Pi_1^{\CA}(\overline{\GA},\overline{\G}, v)=G_{\widehat\CA}$ and so $(\overline{\GA},\overline{\G})$ is an injective profinite graph of pro-$\CA$ groups, since $\overline{\GA}(m)$ are naturally subgroups of $G_{\widehat\CA}$.
		
		To finish the proof we show that the natural map $\alpha': \Gamma\longrightarrow \overline\Gamma$ is an injection. Since $(\GA,\Gamma)$ is reduced $\GA(v)$ and $\GA(w)$ are distinct subgroups in $G$  whenever vertices $v$ and $w$ of $\Gamma$ are distinct. Therefore, as $G$ is residually $\CA$, there exists $U$ such that $v\not\sim_{R_U} w$. This shows that $\alpha'_{|V(\Gamma)}$ is injective. Let $T$ be a maximal subtree of $\Gamma$. Since $T$ has a unique edge connecting two vertices,  $\alpha'_{T}$ is also injective. But the image of the natural map  $\Gamma\longrightarrow \pi_1(\Gamma)$ that sends $T$ to 1   can be viewed as a basis of $\pi_1(\Gamma)$. Since $\pi_1(\Gamma)\leq G$ is residually $\CA$  it follows that $\alpha'_{|E(\Gamma)}$ is injective as well.  
	\end{proof}

	\begin{cor}\label{residually C} Suppose in addition that $G$ contains free  normal subgroup $\Phi$ such that $G/\Phi\in \CA$. Then there exists an open normal subgroup in the pro-$\CA$ topology $V\leq \Phi$ such that for every open normal subgroup $U\leq V$ of $G$ the fundamental group  $\pi_1(\GA_U,\Gamma_U, v_U)$  that appear in (d) is residually $\CA$. In particular,   $(\GA_U, \Gamma_U)$ is injective.  
		
	\end{cor}
	
	\begin{proof} Let $(\overline{\GA},\overline{\G})=\li (\GA_U,\G_U)$, 
		$\Pi_1^{\CA}(\overline{\GA},\overline{\G}, v)=\li (\pi_1(\GA_U,\G_U, v_U))_{\widehat{\CA}}$
		be decompositions of Theorem \ref{thm} (d) and  we can choose all $U\leq \Phi$.  Then the natural epimorphism $\widehat G\longrightarrow G/\Phi$ factors via $   (\pi_1(\GA_V,\G_V, v_V))_{\widehat{\CA}}$ for some $V$ and therefore via every $(\pi_1(\GA_U,\G_U, v_U))_{\widehat{\CA}}$ for $U\leq V$, i.e. we have a homomorphism $\psi_U:(\pi_1(\GA_U,\G_U, v_U))_{\widehat{\CA}}\longrightarrow G/\Phi$ which is injective on vertex groups. Then the kernel of its restriction on $\pi_1(\GA_U,\G_U, v_U)$ is free and so $\pi_1(\GA_U,\G_U, v_U)$ is free-by-$\CA$ group. Thus  $\pi_1(\GA_U,\G_U, v_U)$ is residually $\CA$. 
		
	\end{proof}

	The following theorem covers the remaining part $(d)$ of Theorem \ref{C}.
	
	\begin{theo} \label{std}
		
		We continue with the hypotheses and notation of Theorem \ref{thm}. Furthermore, we assume that $\GA(m)$ is closed in the pro-$\CA$ topology of $G$, for every $m \in \G$. Then the standard tree $S^{abs}=S(\GA,\G)$ of the graph of groups $(\GA,\G)$ is  embedded densely in the standard $\CA$-tree $S=S^{\CA}(\overline{\GA},\overline{\G})$ of the profinite graph of profinite groups $(\overline{\GA},\overline{\G})$ and the action of $\Pi_1^{\CA}(\overline{\GA},\overline{\G})$ on $S^{\CA}(\overline{\GA},\overline{\G})$ extends the natural action of $\pi_1(\GA,\G)$ on $S^{abs}(\GA,\G)$.
	\end{theo}
	
	\begin{proof}
		Let $\widetilde{\G}$, $\widetilde{\overline\Gamma}$ be the abstract universal covering and $\CA$-universal covering of the abstract graph $\G$ and profinite graph $\overline\Gamma$ respectively.  Let $\zeta: \widetilde{\G} \ra \G$, $\overline{\zeta}:\widetilde{\overline{\G}} \ra {\overline{\G}}$   be the covering  and $\CA$-covering maps. Consider  a profinite graph of pro-$\CA$-groups $(\widetilde{\overline{\GA}},\widetilde{\overline{\G}})$ from Definition \ref{universal graph of groups} and its abstract version $(\widetilde{\GA},\widetilde{\G})$.  
		
		If one denotes $\widetilde{G}=\pi_1(\widetilde{\GA},\widetilde{\G},\widetilde{v_0})$, we have the standard tree $\widetilde{S^{abs}}(\widetilde{\GA},\widetilde{\G})$, defined by $\widetilde{S^{abs}}=\bigcup_{\tilde{m} \in \widetilde{\G}} \widetilde{G}/\widetilde{\GA}(\tilde{m})$, $V(\widetilde{S^{abs}})=\bigcup_{\widetilde{v} \in V(\widetilde{\G})} \widetilde{G}/\widetilde{\GA}(\tilde{v})$ and incidence maps $d_0(g\widetilde{\GA}(\tilde{e}))=g\widetilde{\GA}(d_0(\tilde{e}))$ and $d_1(g\widetilde{\GA}(\tilde{e}))=g\widetilde{\GA}(d_1(\tilde{e}))$, $(g \in \widetilde{G}, \tilde{v} \in V(\widetilde{\G}), \tilde{e} \in E(\widetilde{\G}))$.

		In a similar manner, denote $\widetilde\Pi=\Pi_1^{\CA}(\widetilde{\overline \GA}, \widetilde{\overline\Gamma}, \tilde v_0)$ and take its standard $\CA$-tree, that is defined by $\widetilde{S}=\widetilde{S}(\widetilde{\overline{\GA}},\widetilde{\overline{\G}})=\bigcup_{\widetilde{\overline{m}} \in \widetilde{\overline{\G}}} \widetilde{\Pi}/\widetilde{\Pi}(\widetilde{\overline{m}})$, $V(\widetilde{S})=\bigcup_{\widetilde{\overline{v}} \in V(\widetilde{\overline{\G}})} \widetilde{\Pi}/\widetilde\Pi(\widetilde{\overline{v}})$ and incidence maps $d_0(g\widetilde{\Pi}(\widetilde{\overline{e}}))=g\widetilde{\Pi}(d_0(\widetilde{\overline{e}}))$ and $d_1(g\widetilde{\Pi}(\widetilde{\overline{e}}))=g(\widetilde{\Pi}(d_1(\widetilde{\overline{e}}))$, $(g \in \widetilde{\Pi}, \widetilde{\overline{v}} \in V(\widetilde{\overline{\G}}), \widetilde{\overline{e}} \in E(\widetilde{\overline{\G}}))$. 
		
		By Theorem \ref{thm}, $\G$ is densely embedded in $\overline{\G}$ and  $\Pi=G_{\widehat \CA}$ , so there are natural inclusions $\imath_1:\G \ra \overline{\G}$ and $\imath_2: G \ra \Pi$. One defines a morphism of graphs \[\widetilde{\varphi}: \widetilde{S^{abs}}(\widetilde{\GA},\widetilde{\G}) \ra \widetilde{S}(\widetilde{\overline{\GA}},\widetilde{\overline{\G}})\] putting $g\widetilde{\GA}(\tilde{v}) \mapsto g\widetilde{\Pi}(\tilde{v})$, $g\widetilde{\GA}(\tilde{e}) \mapsto g\widetilde{\Pi}(\tilde{e})$, $(g \in \widetilde{G}, \tilde{v} \in V(\widetilde{\G}),\tilde{e} \in E(\widetilde{\G}))$. As $\GA(m)$ $(m \in \G)$ is closed in the pro-$\CA$ topology of $G$ for every $m \in \G$,  $\widetilde{\GA}(m)$ $(\tilde{m} \in \widetilde{\G})$ is also closed in $\widetilde{G}$ for every $\tilde{m} \in \widetilde{\G}$ and the morphism $\widetilde{\varphi}$ is injective. 
		
		Now by \cite[Theorem 6.3.3]{Ribes} there is a canonical isomorphism   $\widetilde{S}\cong S$ and  similarly, $\widetilde{S^{abs}}\cong S^{abs}$ from where we deduce that the $G$-tree $S^{abs}=S(\GA,\G)$  embeds densely in the standard pro-$\CA$ $\Pi$-tree $S=S^{\CA}(\overline{\GA},\overline{\G})$ as needed.
	\end{proof}
	
	Now one uses Theorem \ref{C} (Theorems \ref{thm} and \ref{std}) to prove Theorem \ref{B}:
	
	\begin{proof}[Proof of Theorem \ref{B}]
		As $\pi_1(\GA,\G)$ is residually $\CA$, we can apply Theorem \ref{thm} to conclude part $(a)$. On the other hand, since each $\GA(m)$ is a finite group in $\CA$, it is closed in the pro-$\CA$ topology of $\pi_1(\GA,\G)$, so that the hypotheses of Theorem \ref{std} also hold. This concludes the proof  of $(b)$. 
	\end{proof}
	
	Now we use Theorem \ref{B} to prove the main technical result on interrelation of $S$ and $S^{abs}$ for a virtually free group. We adapt in fact the proof of \cite[Proposition 1.6]{Zalesskii} to the infinitely generated case.
	
	\begin{prop}\label{minimal invariant} Let $G=\pi_1(\GA,\Gamma) $ be the fundamental group of a reduced graph of finite groups   having free subgroup $\Phi$ such that $G/\Phi\in\CA$. Let   $H= \langle h_1, \dots, h_r\rangle$ be an infinite finitely generated   subgroup
		of  $G$  closed in the   pro-$\CA$  topology of $G$, and let $\overline H$ be its closure in the pro-$\CA$ group
		$G_{\widehat \CA}$. Then the standard tree
		$S^{abs}$
		has a unique minimal
		$H$-invariant subtree $D^{abs}$, and its closure $D$ in the standard pro-$\CA$ tree $S$  is the unique minimal $\overline H$-invariant pro-$\CA$
		subtree of $S$; furthermore
		$S^{abs}\cap D= D^{abs}$,  $\overline {D^{abs}}= D$  and $H\backslash D^{abs}= \overline H\backslash D$ is finite.\end{prop}
	
	\smallskip
	\noindent {\it Proof.}   Choose a vertex $v_0$ of $\G$, and denote by $\tilde v_0$ the vertex  $\tilde v_0=
	1\Pi^{abs}(v_0)= 1\Pi (v_0)$  in $S^{abs}\subseteq S$. Define a subgraph $D^{abs}$ of $S^{abs}$ as follows
	$$D^{abs}= \bigcup _{i=1}^rH[\tilde v_0, h_i\tilde v_0].$$
	Put $L= \bigcup _{i=1}^r[\tilde v_0, h_i\tilde v_0]$; this
	is obviously a finite connected graph.  Then  $D^{abs}= HL$. Since $L\cap h_iL\not =\emptyset$  $(i=1,
	\dots, r)$, we have that $D^{abs}$ is a connected subgraph of the tree $S^{abs}$, and so $D^{abs}$ is a tree;
	clearly it is  $H$-invariant. Hence its closure 
	$$D=  \overline {D^{abs}}= \bigcup _{i=1}^r\overline H[\tilde v_0, h_i\tilde v_0] $$ in $S$ is a pro-$\CA$ subtree of $S$; clearly it is $\overline H$-invariant.
	
	Since $H$ is infinite and each $\GA(m)$ is finite, our result will follow by \cite[Lemma 1.4]{Zalesskii} or \cite[Lemma 8.2.1]{Ribes}  after we show that the epimorphism of graphs
	$H\backslash D^{abs}\longrightarrow \overline H\backslash D$ is
	in fact an isomorphism.  To see this we distinguish two cases.
	
	\medskip
	\noindent {\it Case 1.} Assume that $H\le \Phi$.
	
	\smallskip
	\noindent Since the $G$-stabilisers of the elements of $S^{abs}$ are finite groups,  $\Phi$ acts freely on
	$S^{abs}$. By \cite[Lemma 1.1]{Zalesskii} or \cite[Lemma 8.1.1]{Ribes}, there
	exists an open    subgroup   $V$ of $\Phi$ (and so of $G$) containing $H$ such that the map of graphs
	$$H\backslash D^{abs}\longrightarrow H\backslash S^{abs} \longrightarrow V\backslash S^{abs}$$
	is injective. Next
	observe that for every $m\in \G$, we have the following equality of double cosets
	$$V\backslash G/\GA(m) = \overline V\backslash G_{\widehat\CA}/\GA (m)$$
	because  $V$ has finite
	index in $G$; hence, one deduces that
	$V\backslash S^{abs}$ is a subgraph of  $\overline V\backslash S$  from Theorem \ref{B}. Therefore, from the commutative diagram
	$$\xymatrix{H\backslash D^{abs} \ar [r]\ar [d]& H\backslash S^{abs} \ar [d] \ar[r]& V\backslash S^{abs}
		\ar@{^{(}->}[d]\cr
		\overline H\backslash  D\ar[r]&\overline H\backslash   S \ar[r]&\overline V\backslash   S }$$
	we deduce that the left vertical map is injective finishing this case.

	\medskip
	\noindent {\it Case 2.} General case.
	
	\smallskip
	\noindent  Define  $K= \Phi\cap H$. Note that  $K$ is closed  in $\Phi$ and that   $K\backslash D^{abs}$
	is finite (because $K$ has finite index in $H$).  So \cite[Lemma 1.1]{Zalesskii} or \cite[Lemma 8.1.1]{Ribes} can be used again.  Mimicking the argument in Case 1 one shows that $K\backslash D^{abs}=
	\overline K\backslash  D$. What this says is that  if $t, t'\in D^{abs}$, and  $\overline K t= \overline Kt'$, then $K t=
	Kt'$.
	
	Now since  $K$ has finite index in $H$, we have   finite unions  $H= \dbigcup K x_i$ and $\overline H= \dbigcup \overline K x_i$ (for some representatives $x_i\in
	H$ of the cosets $K\backslash H$).   Let  $t, t'\in D^{abs}$, and assume that
	$\overline H t=
	\overline Ht'$. We want to show that then   $H t=  Ht'$.
	By hypothesis we have $\bigcup \overline Kx_it= \bigcup \overline Kx_it'$. So for each $i$, there are some $i'$ and $i''$ such
	that $\overline Kx_it= \overline Kx_{i'}t'$
	and $\overline Kx_it'= \overline Kx_{i''}t$;  hence  by Case 1 $  Kx_it=   Kx_{i'}t'$
	and $  Kx_it'=   Kx_{i''}t$. Therefore, $\bigcup   Kx_it= \bigcup  Kx_it'$, i.e., $Ht=Ht'$.
	\hfill~$\square$
	
	
	\section{Closure of normalisers}
	
	Let $G$ be an abstract group which is residually $\CA$ and let $H$ be a finitely generated closed (in the pro-$\CA$ topology of $G$) subgroup of $G$. In this section we study the relationship between the normaliser $N_G(H)=\{x \in R \mid x^{-1}Hx\}$ of $H$ in $G$ and the normaliser $N_{G_{\widehat{\CA}}}(\overline{H})$ of $\overline{H}$ in $G_{\widehat{\CA}}$. When $G$ contains an open free abstract subgroup, we show (Theorem \ref{A}) that $\overline{N_{G}({H})}=N_{G_{\widehat{\CA}}}(\overline{H})$. This answers Open Question 15.11.10 of \cite{Ribes} and generalises the main result of \cite{Zalesskii}, where the theorem was proved for finitely generated subgroups. In particular, we show that $\overline{N_G(H)}=N_{\widehat G}(\overline H)$ when $G$ is virtually free and $H$ is a finitely generated subgroup of $G$ (Corollary \ref{cor2}).

	We state Theorem \ref{A} in a more general form, namely for the pro-$\CA$ case.
	
	\begin{theo} \label{a}
		Let $G$ be a group having normal free subgroup $\Phi$ such that $G/\Phi\in \CA$. Let $H$ be a finitely generated subgroup of $G$.  Then   $\overline{N_G(H)}=N_{ G_{\widehat\CA}}(\overline{H})$, where the closure $\overline{N_G(H)}$ is taken in $G_{\widehat\CA}$.
	\end{theo}
	
	\begin{proof}
		Obviously $\overline{N_G(H)} \leq N_{G_{\widehat{\CA}}}(\overline{H}).$ We need to prove the opposite containment. 
		
		We continue with the notations of Sections 2 and 3. According to a result of Scott (cf. \cite{Scott}), $G$ is the abstract fundamental group $\pi_1(\GA,\G)$ of a graph of groups $(\GA,\G)$ over a graph $\G$ such that each $\GA(m)$ $(m \in \G)$ is a finite group. A finite subgroup of $G$ is isomorphic to a subgroup of $G/\Phi$, and so $\GA(m)$ is in $\CA$.  By Theorem \ref{thm} (d) the graph of groups $(\overline{\GA},\overline{\G})$ decomposes as a surjective inverse limit $(\overline{\GA},\overline{\G})=\li (\GA_U,\G_U)$ of finite graphs of finite $\CA$-groups and  
		$
		\Pi_1^{\CA}(\overline{\GA},\overline{\G}, v)=\li (\pi_1(\GA_U,\G_U, v_U))_{\widehat{\CA}}$. 
		Let $\pi_U:\Pi_1^{\CA}(\overline{\GA},\overline{\G}, v)\longrightarrow  (\pi_1(\GA_U,\G_U, v_U))_{\widehat{\CA}}$ be the natural projection. Since $\pi_U$ induced by the morphism $(\overline{\GA},\overline{\G}, v)\longrightarrow (\GA_U,\G_U, v_U)$, one has $\pi_U(\pi_1(\GA,\Gamma))\leq \pi_1(\GA_U, \Gamma_U)$. Put $G_U=\pi_1(\GA_U, \Gamma_U, v_U)$ and  $H_U=Cl(\pi_U(H))$, where $Cl$ means the closure in the pro-$\CA$ topology of $G_U$, and note that by Corollary \ref{residually C} $G_U$ is residually $\CA$.
		Then by \cite[Theorem 13.1.7]{Ribes} $\overline{N_{G_U}(H_U)}=N_{ (G_U)_{\widehat\CA}}(\overline{H_U})$. Since $\overline H_U=\overline{\pi_U(H)}$, $N_{G_{\widehat\CA}}(\overline{H})=\li_U N_{ (G_U)_{\widehat\CA}}(\overline{H_U})$ and $\overline{N_G(H)}=\li_U \overline{N_{G_U}(H_U)}$ one deduces that $\overline{N_G(H)}=N_{ G_{\widehat\CA}}(\overline{H})$.
	\end{proof}
	
	In \cite{Hall}, Marshall Hall proved that a finitely generated subgroup $H$ of a free abstract group $\Phi$ is closed in the profinite topology of $\Phi$. It easily follows that a finitely generated subgroup of a virtually free (or free-by-finite) abstract group $G$ is automatically closed in the profinite topology of $G$. Therefore one has the following corollary of Theorem \ref{A}.
	
	\begin{cor} \label{cor2}
		Let $G$ be a virtually free abstract group and let $H$ be a finitely generated subgroup of $G$. Then
		$$\overline{N_G(H)}=N_{\hat{G}}(\overline{H}).$$
	\end{cor}

	We shall finish the section proving the same equality for centralizers of cyclic subgroups that generalizes \cite[Corollaries 13.10, 13.1.12]{Ribes} . We need however first the following 	
	
	\begin{lemma} \label{cor}
		Let $G$ be a group having a normal free subgroup $\Phi$ of $G$ such that $G/\Phi \in \CA$. Let $H$ be a cyclic non-trivial subgroup of $\Phi$. Then $$C_G(H)=C_G(Cl(H)),$$ where $Cl(H)$ denotes the closure of $H$ in the pro-$\CA$ topology of $G$.
	\end{lemma}
	
	\begin{proof}
		As $\Phi$ is closed in $G$, $Cl(H)$ is also the closure of $H$ in the pro-$\CA$ topology of $\Phi$. By \cite[Proposition 3.4]{RZ94}, $Cl(H)$ is cyclic and contains $H$ as a subgroup of finite index. Say $Cl(H)=\langle x \rangle$ and $H=\langle x^n \rangle$. Now, if $a \in G$ and $a^{-1}x^na=x^n$, then both $a^{-1}xa$ and x are n-th roots of $x^n$. Since in a free abstract group n-th roots are unique, we deduce that $a^{-1}xa=x$ and the result follows.
	\end{proof}
	
	\begin{prop}
		Let $G$ be a free-by-$\CA$ abstract group. If $H$ is an infinite cyclic subgroup of $G$, then $$C_{G_{\hat{C}}}(\overline{H})=\overline{C_{G}(H)}.$$ 
	\end{prop}
	
	\begin{proof}
		Consider the natural homomorphism $$\varphi:N_G(H) \ra Aut(H) \cong \mathbb{Z}/2\mathbb{Z}.$$ Then $Ker(\varphi)=C_G(H)$.
		
		Assume first that $H$ is closed. Note that $$\overline{C_G(H)} \leq C_{G_{\widehat{\CA}}} \leq N_{G_{\widehat{\CA}}}(\overline{H})=\overline{N_G(H)}$$ (for the last equality we use Theorem \ref{A}). Since the index of $\overline{C_G(H)}$ in $\overline{N_G(H)}$ is at most $2$, the result follows easily: suppose $C_{G_{\widehat{\CA}}}(\overline{H})=\overline{N_G(H)}$ and let $g \in N_G(H)$; then $g \in C_{G_{\widehat{\CA}}}(\overline{H})$, and so $g \in C_G(H)$, i.e., $C_G(H)=N_G(H)$. Hence $\overline{C_G(H)}=C_{G_{\widehat{\CA}}}(\overline{H})$.
		
		Assume now that $H$ is a cyclic subgroup of $\Phi$, not necessarily closed. By Lemma \ref{cor}, $C_G(H)=C_G(Cl(H))$. Therefore using the result above for the closed subgroup $Cl(H)$, $$\overline{C_G(H)}=\overline{C_G(Cl(H))}=C_{G_{\widehat{\CA}}}(\overline{Cl(H)})=C_{G_{\widehat{\CA}}}(\overline{H}),$$ since $\overline{H}=\overline{Cl(H)}$.
	\end{proof}

	\section{Subgroup conjugacy separability}
	
	A group $G$ is said to be subgroup conjugacy $\CA$-separable if whenever $H_1$ and $H_2$ are finitely generated closed subgroups of $G$ (in its pro-$\CA$ topology), then $H_1$ and $H_2$ are conjugate in $G$ if and only if their images in every quotient $G/N \in \CA$ are conjugate, or equivalently for residually $\CA$ groups, if and only if their closures in $G_{\widehat{\CA}}$ are conjugate. In this section we prove the subgroup conjugacy $\CA$-separability of a  free-by-$\CA$  group $G$ (Theorem \ref{D}). This answers Open Question 15.11.11 of \cite{Ribes} and generalises the main result of \cite{Chagas} where it is proved for finitely generated free-by-$\CA$  groups. In particular, if $G$ is virtually free then $G$ is subgroup conjugacy separable (Corollary \ref{cor2}). One uses again the technique  developed in Section 3 of groups acting on trees and the interrelation between abstract and profinite graphs and groups.
	
	Thus for the rest of the section fix a group $G$ having a free subgroup $\Phi$ with $G/\Phi\in \CA$. 	According to \cite{Scott},  $G$ splits as the fundamental group of a graph of groups $(\GA,\G)$ over a graph $\G$, i.e., $G=\pi_1(\GA,\G)$, such that $\GA(m) \in \CA$ for every $m \in \G$. In fact we may assume that $(\GA,\G)$ is reduced, i.e. whenever $e$ is an edge of $\G$ which is not a loop, then the order of the finite group $\GA(e)$   is strictly smaller than the order of $\GA(d_i(e))$ $(i=0,1)$.  Then  $G_{\widehat{\CA}}$ is the pro-$\CA$ fundamental group of the profinite graph of groups $(\overline{\GA},\overline{\G})$ constructed in Theorem \ref{thm}.
	
	\begin{lemma} \label{rgg2}
		$(\overline{\GA},\overline{\G})$ is reduced, i.e. $\GA(e)\neq \GA(d_i(e))$ for all $e\in \overline\Gamma$, $i=1,2$.
	\end{lemma}
	
	\begin{proof}
		As  in the proof of Theorem \ref{C}, let $R_U$ be the equivalence relation in $\G$ defined by $v,w \in V(\G)$, $v \sim_{R_U} w$ if $\GA_U(v)=\GA_U(w)$ and given $e,e' \in E(\G)$, $e \sim_{R_U} e'$ if $d_0(e) \sim_{R_U} d_0(e')$ and $d_1(e) \sim_{R_U} d_1(e')$, where all $U$ can be taken inside $\Phi$ and $\GA_U(m)=\GA(m)U/U$. Hence, the quotient graph $\G_U$ defined by $\G_U=\G/R_U$ is finite and we have the finite graph of groups $(\GA_U,\G_U)$ such that $\GA_U(\bar m)=\GA(m)$ where $\bar m$ is the equivalence class of $m$ and $\partial_i(gU/U)=\partial_i(g)U/U$.
		
		Therefore, if for given $e \in E(\G)$ which is not a loop  neither $\partial_1:\GA(e) \ra \GA(d_1(e))$ nor $\partial_0:\GA(e) \ra \GA(d_0(e))$ is an isomorphism, then neither $\partial_1:\GA_U(e) \ra \GA_U(d_1(e))$ nor $\partial_0:\GA_U(e) \ra \GA_U(d_0(e))$ can be an isomorphism, so each finite graph of groups $(\GA_U,\G_U)$ is reduced. Moreover, for a   subgroup $V\leq U$ of finite index normal in $G$ the morphism $\eta_{VU}:(\GA_V,\G_V)\longrightarrow (\GA_U,\G_U)$ restricted on each vertex group  $\GA(v)$ is an isomorphism.   Hence the inverse limit preserve the  property of being reduced and the proof is finished.  
	\end{proof}
	
	\begin{proof}[Proof of Theorem \ref{D}]
		We continue with the notation of the section.	Let $H_1$ and $H_2$ be finitely generated closed subgroups of $G$. Since $G$ is residually $\CA$, it suffices to prove that if $\gamma \in G_{\widehat{\CA}}$ and $\overline{H_1}=\gamma\overline{H_2}\gamma^{-1}$, then there exists some $g \in G$ such that $H_1=gH_2g^{-1}$. 
		
		As usual, we denote by $S^{abs}$ and $S$ the abstract standard tree of $(\GA,\G)$ and the standard pro-$\CA$-tree of $(\overline{\GA},\overline{\G})$.
		We divide the proof in two cases, when $H_1$ is infinite and when $H_1$ is finite.
		
		\medskip
		Case $1.$ $H_1$ is infinite (hence so is $H_2$). By Proposition \ref{minimal invariant}, $S^{abs}$ has a unique minimal $H_i$-invariant subtree $D_i^{abs}$, and $D_i=\overline{D_i^{abs}}$ is the unique minimal $\overline{H_i}$-invariant pro-$\CA$ tree of $S$ $(i=1,2)$. Then $\gamma D_2$ is a minimal $\overline{H_1}$-invariant pro-$\CA$ tree of $S$, and hence $D_1=\gamma D_2$.
		
		By Theorem \ref{std} and Proposition \ref{minimal invariant} one has that $D_i^{abs}$ is a connected component of $D_i$ considered as an abstract graph and any other component of $D_i$ has the form $\beta D_i^{abs}$, for some $\beta \in \overline{H_i}$ $(i=1,2)$. Therefore $\gamma D_2^{abs}$ is an abstract connected component of $D_1$. It follows that there exists some $\tilde{h}_1 \in \overline{H_1}$ such that \[\tilde{h}_1\gamma D_2^{abs}=D_1^{abs}.\] Since $H_2$ is infinite and the $G$-stabiliser of any $m \in S^{abs}$ is finite, the tree $D_2^{abs}$ must contain at least one edge; say $e \in E(D_2^{abs})$. Then $\tilde{h}_1\gamma e \in D_1^{abs} \subseteq S^{abs}$. Since $G \bs S^{abs}=\G\subseteq \overline \Gamma$ by Theorem \ref{B}, there exists some $g_1 \in G$ such that $g_1e=\tilde{h}_1\gamma e$. 
		Hence $g_1^{-1}\tilde{h}_1\gamma$ is in the $G_{\widehat{\CA}}$-stabiliser of $e$, which in fact coincides with the $G$-stabiliser of $e$ since it is finite. 
		Therefore $g_1^{-1}\tilde{h}_1\gamma \in G$, and so $\tilde{h}_1\gamma=g \in G$. 
		Finally, taking into account that $H_1$ and $H_2$ are closed, we have 
		\begin{equation*}
			H_2=G \cap \overline{H_2}=G \cap \overline{H_1}^{\gamma}=G \cap\overline{H_1}^{(\tilde{h}_1)^{-1}g}=G \cap \overline{H_1^{g}}=H_1^g,
		\end{equation*}	
		as desired. 
		
		Case $2.$ $H_1$ is finite (hence so is $H_2$).  Since $\overline{\Phi}$ is open in $G_{\widehat{\CA}}$,  $G_{\widehat{\CA}}=G\overline{\Phi}$. Hence $\gamma=g'\gamma'$, where $g' \in G$ and $\gamma' \in \overline{\Phi}$. Then $H_2=H_1^{\gamma}=(H_1^{g'})^{\gamma'}$. So, replacing $H_1$ with $H_1^{g'}$ and $\gamma$ with $\gamma'$, we may assume that $\gamma \in \overline{\Phi}$. Hence $\overline{\Phi}H_1=\overline{\Phi}H_2$, and so $\Phi H_1=\Phi H_2$, because $\Phi$ is closed in $G$. Since in fact $\Phi H_i$ is open in $G$, one has $(\Phi H_i)_{\widehat{\CA}} \leq G_{\widehat{\CA}}$ $(i=1,2)$. Thus from now on we may assume that \[G=\Phi H_1=\Phi H_2=\Phi \rtimes H_1 = \Phi \rtimes H_2.\] This implies that $H_1$ and $H_2$ are maximal finite subgroups of $G$, and so they are $G$-stabilisers of some vertices $S^{abs} \subseteq S$, say $v_1$ and $v_2$, respectively (see, for example, \cite[Theorem I.15]{Serre} and \cite[Theorem 7.1.2]{Ribes} or \cite[Theorem 2.10]{Melnikov}). 
		
		Recall that by Lemma \ref{rgg2}, $(\overline{\GA},\overline{\G})$ is reduced.
		Then, if $v$ is one of the vertices of $\tilde{e}$ and  $\overline{G}_v=\overline{G}_{\tilde{e}}$, it must be because the image of $\tilde e$ in $\overline\Gamma$ is a loop.
		
		Since $H_1$ stabilises $\gamma v_2$ and $v_1$, it must stabilise every element of the chain $[\gamma v_2,v_1]$ in $S$ (see \cite[Corollary 4.1.6]{Ribes}); therefore, since $H_1$ is maximal, it is the $\overline{G}$-stabiliser of each element of the chain $[\gamma v_2,v_1]$. In particular, the end points of any edge of this chain have the same $\overline{G}$-stabilisers. By the comment above, this means that the projection on $\overline\G$ of any edge of $[\gamma v_2,v_1]$ must be a loop. Since the image of $[\gamma v_2,v_1]$ in $\overline\G$ is a connected subgraph of the graph $\overline\G$, one deduces that the image of $[\gamma v_2,v_1]$ in $\overline\G$ has a unique vertex. Therefore $\gamma v_2$ and $v_1$ are in the same $\overline{G}$-orbit. Hence $\overline{G}v_2=\overline{G}\gamma v_2=\overline{G} v_1$.
		
		Since $G \bs S^{abs}=\G$ is densely embedded in $\overline{\G}=\overline{G} \bs S$, one deduces that $G v_2=G v_1$. Say $gv_2=v_1$, where $g \in G$. Then $H_1=G_{v_1}=gG_{v_2}g^{-1}=gH_2g^{-1}$.  
	\end{proof}
	
	By the  M. Hall theorem every finitely generated subgroup of a virtually free group is closed in the profinite topology. Thus we have the following 
	
	\begin{cor}
		Let $G$ be a virtually free group. Then $G$ is subgroup conjugacy separable.
	\end{cor}

	\bigskip
	{\it Author's Adresses:}
	
	\medskip
	Mattheus P. S. Aguiar\\
	Departamento de Matem\'atica,\\
	~Universidade de Bras\'\i lia,\\
	70910-900 Bras\'\i lia DF,\\
	Brazil
	
	mattheus@mat.unb.br
	
	\medskip
	Pavel A. Zalesski\\
	Departamento de Matem\'atica,\\
	~Universidade de Bras\'\i lia,\\
	70910-900 Bras\'\i lia DF\\
	Brazil
	
	pz@mat.unb.br

\end{document}